\documentclass[a4paper,reqno,10pt]{amsart}
\usepackage[colorlinks=true,linkcolor=blue,citecolor=blue,urlcolor=blue]{hyperref}
\usepackage{a4wide,amssymb,amsmath,amsthm,amsfonts,mathtools,tikz-cd,enumitem}
\newtheorem{theorem}{Theorem}[section]
\newtheorem{proposition}[theorem]{Proposition}
\newtheorem{lemma}[theorem]{Lemma}
\newtheorem{corollary}[theorem]{Corollary}
\newtheorem*{theoremA}{Theorem A}
\theoremstyle{definition}
\newtheorem{definition}[theorem]{Definition}
\newtheorem{construction}[theorem]{Construction}
\newtheorem{setting}[theorem]{Setting}
\newtheorem{remark}[theorem]{Remark}
\newtheorem{example}[theorem]{Example}
\newtheorem*{organization}{Organization}
\newtheorem*{conventions}{Conventions and notation}
\newtheorem*{formalization}{Formalization}
\newtheorem*{aiuse}{Use of AI}
\DeclareMathOperator{\Hom}{Hom}
\DeclareMathOperator{\End}{End}
\DeclareMathOperator{\Ext}{Ext}
\DeclareMathOperator{\rad}{rad}
\DeclareMathOperator{\im}{im}
\DeclareMathOperator{\supp}{supp}
\DeclareMathOperator{\Mod}{Mod}
\DeclareMathOperator{\md}{mod}
\newcommand{\dm}{\dim_k}
\newcommand{\D}{\mathcal D}
\newcommand{\C}{\mathcal C}
\numberwithin{equation}{section}
\setlist[enumerate,1]{label=\textup{(\arabic*)},ref=\arabic*}
\title[Representation embeddings and the third Brauer--Thrall conjecture]
{Representation embeddings and the third Brauer--Thrall conjecture}
\author[H.\ Enomoto]{Haruhisa Enomoto}
\address{Parakeet Inc., Japan}
\email{haruhisa.enomoto.math@gmail.com}
\date{\today}
\begin{document}
\begin{abstract}
Let $A$ be a representation-infinite finite-dimensional algebra over an
algebraically closed field $k$. We prove that for every infinite cardinal
$\kappa$ there is an indecomposable $A$-module $M$ with $\dim_k M=\kappa$.
Over algebraically closed fields, this settles Simson's conjecture, also
called the third Brauer--Thrall conjecture, that every
representation-infinite finite-dimensional algebra has indecomposable
modules of arbitrarily large cardinality. To construct these modules, we
extend Bongartz's representation embedding from finite-dimensional
Kronecker modules without a simple injective summand to such modules of
arbitrary dimension. The extended
functor commutes with direct sums and filtered colimits and preserves the
dimension of infinite-dimensional modules. Applying it to the
indecomposable Kronecker modules of every infinite dimension constructed by
G\"obel and Simson yields the required $A$-modules.
\end{abstract}
\maketitle
\raggedbottom
\begingroup\small\tableofcontents\endgroup

\section{Introduction}
Let $A$ be a finite-dimensional algebra over a field $k$. The algebra $A$
is called \emph{representation-finite} if there are only finitely many
isomorphism classes of finite-dimensional indecomposable $A$-modules, and
\emph{representation-infinite} otherwise. The Brauer--Thrall conjectures
concern the indecomposable modules of a representation-infinite
algebra; see Ringel \cite{RingelReport} for their origin. The first
Brauer--Thrall conjecture asserts that a representation-infinite algebra
has finite-dimensional indecomposable modules of arbitrarily large
dimension. Roiter \cite{Roiter} proved it for finite-dimensional algebras
over any field, and Auslander \cite{AuslanderII} proved the analogue for
Artin rings, with length in place of dimension. The second Brauer--Thrall conjecture asserts that,
over an infinite field, a representation-infinite algebra has infinitely
many positive integers $d$ for which there are infinitely many pairwise
non-isomorphic indecomposable modules of dimension $d$. Over algebraically
closed fields of arbitrary characteristic, it was proved by Bautista \cite{BautistaBT2}, using
Bongartz's standardness theorem \cite{BongartzStandard};
Fischbacher \cite{Fischbacher} and Bretscher and Todorov
\cite{BretscherTodorov} subsequently gave alternative proofs (see
Bongartz \cite[\S5.1]{Bongartz} for the history). In \cite[Theorem~6]{Bongartz}, Bongartz also proved the second conjecture
over algebraically closed fields by constructing representation embeddings
from the category of finite-dimensional Kronecker modules
without a simple injective summand into the category $\md A$ of finite-dimensional $A$-modules. These functors transfer indecomposable
Kronecker modules to indecomposable $A$-modules and are the starting
point of this paper.

Both of these conjectures concern finite-dimensional modules. For modules of arbitrary size,
Ringel and Tachikawa \cite[Corollary~4.4]{RingelTachikawa} proved that every module over
a representation-finite Artin algebra is a direct sum of finitely
generated indecomposable modules (see also Auslander \cite[pp.~1--2]{AuslanderLarge}), and
Auslander \cite[Theorem~A]{AuslanderLarge} proved conversely that if every
indecomposable module over an Artin algebra is finitely generated, then
the algebra is representation-finite. Thus every representation-infinite
finite-dimensional algebra has infinite-dimensional indecomposable
modules, and one may ask how large they can be.
As an infinite-dimensional counterpart of the first Brauer--Thrall conjecture,
Simson \cite[Conjecture~$1_\infty$]{SimsonCorner}
conjectured that
every representation-infinite finite-dimensional algebra over a field has
indecomposable modules of arbitrarily large cardinality, meaning that for every
infinite cardinal $\kappa$ there exists an indecomposable module $M$
with $|M|\ge\kappa$. Following Iovanov \cite[\S6.1]{Iovanov}, who calls this a
Brauer--Thrall~3 question, we refer to it as the \emph{third Brauer--Thrall
conjecture}. Since any infinite-dimensional
vector space $M$ over $k$ satisfies $|M|=\max(|k|,\dm M)$, the condition
that $|M|$ can be made arbitrarily large is equivalent to the existence, for
every infinite cardinal $\kappa$, of an indecomposable $A$-module $M$ with
$\dm M\ge\kappa$.

In this paper, we settle this conjecture over algebraically closed fields by proving that every representation-infinite
finite-dimensional algebra over an algebraically closed field has
an indecomposable module of dimension equal to each prescribed infinite
cardinal. Our main result (Theorem~A) constructs representation embeddings on the
subcategory of Kronecker modules of arbitrary dimension without a simple injective summand. Applying these
embeddings to the indecomposable Kronecker modules of every infinite dimension
constructed by G\"obel and Simson \cite{GoebelSimson} gives
Corollary~\ref{cor:dimension}, which also counts the resulting modules and
bounds the nilpotence exponent of their endomorphism radicals.

\subsection{Representation embeddings on modules of arbitrary dimension}
Let $K_2$ be the Kronecker quiver $1\rightrightarrows2$.
We consider its representations of arbitrary dimension, identified with
left $kK_2$-modules. Write such a representation as $(V,W;a,b)$, with $V$
at vertex~$1$, $W$ at vertex~$2$ and arrow maps $a,b:V\to W$. Let $S_+=(k,0;0,0)$ and $S_-=(0,k;0,0)$ be its simple
injective and simple projective representations, respectively.
Let $\D$ be the full subcategory of $\Mod kK_2$ consisting of the modules
$T$ without a direct summand isomorphic to $S_+$, and put $\D_{\mathrm{fin}}=\D\cap\md kK_2$.
Since $S_+$ is simple
and injective, $T\in\D$ if and only if $\Hom_{K_2}(S_+,T)=0$; the category $\D$
is closed under submodules, direct sums and filtered colimits
(Lemma~\ref{lem:Kronecker}(2)), and a short sequence in $\D$ is
called exact when it is short exact in $\Mod kK_2$.
The finite-dimensional subcategory $\D_{\mathrm{fin}}$ is the common domain of Bongartz's embeddings:
one of his constructions, applied to all Kronecker modules, sends $S_+$ to zero and therefore cannot reflect
isomorphism classes on that larger category.

Following Bongartz \cite[\S2.1]{Bongartz}, we call a
$k$-linear functor $G:\D_{\mathrm{fin}}\to\md A$ or $G:\D\to\Mod A$ a \emph{representation embedding} if it preserves short exact sequences and
indecomposability and reflects isomorphism classes, in the sense that $G(T)\cong G(U)$
implies $T\cong U$.
Our main result is the following.
\begin{theoremA}
Let $A$ be a representation-infinite finite-dimensional algebra over an
algebraically closed field $k$. There is a representation embedding
\[
 F:\D\longrightarrow\Mod A
\]
with the following properties.
\begin{enumerate}
\item $F$ sends finite-dimensional modules to finite-dimensional modules.
\item $F$ commutes with arbitrary direct sums and filtered colimits.
\item If $T\in\D$ has infinite dimension, then $\dm F(T)=\dm T$.
\end{enumerate}
\end{theoremA}
Theorem~\ref{thm:embedding} proves a stronger form for all the choices of
Bongartz's construction described in Section~\ref{sec:outline}, and also
shows that the resulting functors are faithful.

\subsection{Indecomposable modules of prescribed infinite dimension}
For every infinite cardinal $\kappa$, G\"obel and Simson
\cite[Theorem~1.2]{GoebelSimson} construct a family of Kronecker modules $T_\beta$
($\beta\subseteq\kappa$) of dimension $\kappa$, with neither simple module as a direct summand, satisfying
\[
 \Hom_{K_2}(T_\beta,T_\gamma)\cong
 \begin{cases}k,&\beta\subseteq\gamma,\\0,&\beta\nsubseteq\gamma,\end{cases}
 \qquad \End_{K_2}(T_\beta)\cong k
\]
(see Theorem~\ref{thm:rigid-family}).
For our embeddings $F$, we construct algebra surjections
$\End_A(F(T_\beta))\to\End_{K_2}(T_\beta)\cong k$ that are left inverses
to the maps induced by $F$. Their kernels are nilpotent, so the endomorphism
algebras $\End_A(F(T_\beta))$ are local with residue field $k$.
Using Bongartz's bound \cite[Theorem~7]{Bongartz} on the
rank over $k[t]$ of a free bimodule defining an embedding of
finite-dimensional $k[t]$-modules, we can choose $F$ so that the nilpotence
exponent of these kernels is at most
$\max(30,4\dm A)$ for every $\kappa$ and every $\beta\subseteq\kappa$
(Proposition~\ref{prop:global-refinement}). Combining Theorem~A with the
family of G\"obel and Simson then yields the following result.
\begin{corollary}[= Corollary~\ref{cor:radical}]\label{cor:dimension}
Let $A$ be a representation-infinite finite-dimensional algebra over an
algebraically closed field $k$. There is an integer
$1\le n\le\max(30,4\dm A)$ such that, for every infinite cardinal
$\kappa$, there are $2^\kappa$ pairwise non-isomorphic indecomposable
$A$-modules $M$ satisfying
\[
 \dm M=\kappa,\qquad
 \End_A(M)/\rad\End_A(M)\cong k,\qquad
 (\rad\End_A(M))^n=0.
\]
If $\kappa\ge |k|$, this is the largest possible number of isomorphism
classes in dimension $\kappa$.
\end{corollary}

\subsection{Outline of the proof}
\label{sec:outline}
We first reduce to a basic minimal representation-infinite algebra,
that is, a basic representation-infinite algebra whose proper quotients
are representation-finite. For such an algebra, Bongartz \cite[Theorem~6 and its proof]{Bongartz} supplies at least
one of three constructions of a representation embedding on $\D_{\mathrm{fin}}$ (Theorem~\ref{thm:finite-input}):
a quotient of projective modules, a push-down along a Galois covering, or a
push-down along a closed walk. Because these constructions are defined using tensor products with
finite-dimensional spaces, quotients, and finite direct sums of vertex spaces,
their formulas extend directly to modules of arbitrary dimension; we establish their exactness and colimit properties
in Section~\ref{sec:constructions}.

To control homomorphisms between the resulting modules, let $\C$ be the
full subcategory of $\D$ consisting of the Kronecker modules with
neither simple module as a direct summand. For each of the three functors
$F$, we construct $k$-linear maps
\[
 \rho_{T,U}:\Hom_A(F(T),F(U))\longrightarrow\Hom_{K_2}(T,U),
 \qquad T,U\in\C,
\]
preserving composition and satisfying $\rho_{T,U}(F(f))=f$, so that each $\rho_{T,T}$ is a surjective algebra homomorphism onto $\End_{K_2}(T)$.
For each $T$, the kernel of $\rho_{T,T}$ is a nilpotent ideal of $\End_A(F(T))$,
with exponent bounded independently of $T$. Consequently, when $T\in\C$ is indecomposable,
$\End_A(F(T))$ contains no nontrivial idempotent and $F(T)$ is indecomposable;
likewise, applying $\rho_{T,U}$ and $\rho_{U,T}$ to an isomorphism $F(T)\xrightarrow{\sim}F(U)$ and its inverse recovers an
isomorphism $T\xrightarrow{\sim}U$. Finally, to extend reflection of isomorphism classes from $\C$ to all
of $\D$, we split off the simple projective summands $S_-$ and cancel their images under $F$.

\subsection{Related work}
\label{sec:intro-related}
Classical constructions of infinite-dimensional indecomposable modules include generic,
Pr\"ufer and adic modules; see Crawley-Boevey
\cite[\S\S1.3--1.5, 3.1]{CrawleyBoeveyInfinite}. These modules
do not by themselves attain every infinite dimension; for instance,
the Pr\"ufer modules over the Kronecker algebra all have countable
dimension.
For a fully $k$-wild algebra $A$, meaning that there exists a fully
faithful exact $k$-linear functor from finite-dimensional
$k\langle x,y\rangle$-modules to finite-dimensional $A$-modules,
Simson \cite[Proposition~2.2]{SimsonCorner} constructs a fully faithful
exact $k$-linear functor $\Mod k\langle x,y\rangle\to\Mod A$ and obtains
indecomposable modules of arbitrarily large infinite dimension.
Over an algebraically closed field, Simson \cite[Theorem~2.8]{SimsonLarge} also constructs indecomposable modules of
arbitrarily large infinite dimension for four classes of
representation-infinite algebras:
loop-finite algebras, strongly simply connected algebras, algebras with
radical square zero, and group algebras of finite groups.
By contrast, Corollary~\ref{cor:dimension} provides indecomposable modules of every infinite dimension over any
representation-infinite finite-dimensional algebra over an algebraically closed field.

In \cite[Theorem~3.1]{SimsonCorner}, Simson states that every
representation-infinite special biserial algebra over an arbitrary field $k$
has indecomposable modules of arbitrarily large cardinality.
However, Case~4 of that proof relies on \cite[Proposition~3.7]{SimsonCorner}, which asserts that
every such algebra $R$ has a corner $eRe$ with a quotient
isomorphic either to $k\langle x,y\rangle/(x^2,xy,yx,y^2)$ or to the
path algebra of an acyclic orientation of a cycle with at least two vertices.
Example~\ref{ex:simson-corner} in Appendix~\ref{sec:simson-counterexample}
refutes this assertion over every field, leaving a gap in Case~4.
Over algebraically closed fields, Corollary~\ref{cor:dimension} establishes
the existence of indecomposable modules of arbitrarily large cardinality for all representation-infinite special biserial algebras; over other fields, that conclusion
is not established by Simson's proof.

For a representation-infinite finite-dimensional algebra $R$ over a
field $K$ with $(\rad R)^2=0$ and $R/\rad R\cong K^m$ for some $m\ge1$,
Simson \cite[Corollary~3.6]{SimsonCorner} constructs $R$-modules $U$ of arbitrarily large cardinality satisfying
\[
 \End_R(U)/\rad\End_R(U)\cong K,\qquad
 (\rad\End_R(U))^2=0,
\]
so that $\End_R(U)$ is local, although the dimension of $U$ is not
prescribed.
Over an algebraically closed field, Corollary~\ref{cor:dimension} bounds the nilpotence exponent of the endomorphism radical uniformly
for families of indecomposable $A$-modules in every infinite dimension.

Ringel \cite[Corollary on p.~407]{RingelInfinite} showed that for a
connected tame hereditary algebra $R$, the category of all $R$-modules
contains a full exact subcategory equivalent to the category of all
modules over a wild hereditary algebra; thus tame hereditary
algebras are wild with respect to modules of arbitrary dimension. For
representation-infinite algebras $A$ over an algebraically closed field
$k$ that are loop-finite or
strongly simply connected, Simson \cite[Theorem~1.2]{SimsonRealization}
constructs a fully faithful exact functor from the category of all
Kronecker modules to $\Mod A$ sending finite-dimensional modules to
finite-dimensional modules. In contrast to these full embeddings for
special classes, Theorem~A applies to every representation-infinite algebra
over an algebraically closed field, but its functor on $\D$ need not be full,
as Example~\ref{ex:not-full} shows.

Iovanov \cite[\S6.1]{Iovanov} asks whether a given exact functor between
categories of all comodules whose restriction to finite-dimensional
comodules is a representation embedding must itself be a representation
embedding. For modules over finite-dimensional algebras,
Theorem~\ref{thm:embedding} establishes a version with domain $\D$:
for the functors constructed here, being a representation embedding on
$\D_{\mathrm{fin}}$ implies being a representation embedding on $\D$.

Gregory and Prest \cite[Theorem~4.2]{GregoryPrest} study tensor representation embeddings
$G:\Mod S\to\Mod R$ induced by finite-dimensional bimodules between
finite-dimensional algebras, and construct a functor $H:\Mod R\to\Mod S$ satisfying
$H(G(X))\cong X$ for every $S$-module $X$.
Their construction assumes that for finite-dimensional $S$-modules $X$ and $Y$, the space $\Hom_R(G(X),G(Y))$
decomposes as the direct sum of the image of $\Hom_S(X,Y)$ and an ideal of radical morphisms factoring through
finite direct sums of a fixed finite-dimensional $R$-module. By contrast, our maps $\rho_{T,U}$ are defined only between objects in the image
of the embedding rather than as a functor on all of $\Mod A$, and in place
of this factorization condition we require only a uniform nilpotence bound on
the endomorphism kernels $\ker\rho_{T,T}$
(Definition~\ref{def:retraction}).

Over perfect fields, Bautista, P\'erez and Salmer\'on
\cite[Theorem~1.4]{BautistaPerezSalmeron} construct a representation
embedding into $\md A$ for every representation-infinite
finite-dimensional algebra $A$, with domain the category of modules over a
principal ideal domain that have finite dimension over $k$.
While this provides a finite-dimensional representation embedding over a wider class of fields,
our argument relies on the specific quotient, covering and walk constructions
of Theorem~\ref{thm:finite-input}, which are established over
algebraically closed fields.

\begin{organization}
In Section~\ref{sec:extension}, we describe how to split off the simple
summands of a Kronecker module and show that retractions on Hom spaces with nilpotent
kernels preserve indecomposability and reflect isomorphisms. In Section~\ref{sec:constructions}, we
reduce to basic minimal representation-infinite algebras, define the three functors in arbitrary dimension, and state
Theorem~\ref{thm:finite-input}. In Section~\ref{sec:control}, we construct
retractions on Hom spaces and bound the nilpotence exponents of their kernels using the Loewy length in the quotient
case and the fibre sizes in the covering and walk cases.
In Section~\ref{sec:completion}, we cancel the simple projective summands, transfer the
embedding to the original algebra, and prove the uniform numerical bound and
Corollary~\ref{cor:dimension}. Appendix~\ref{sec:simson-counterexample} gives
the counterexample discussed in Section~\ref{sec:intro-related}.
\end{organization}
\begin{conventions}
Throughout, the letter $k$ denotes a field, assumed to be algebraically closed where stated, and we assume the axiom of choice.
All algebras and categories are $k$-linear, and modules are left modules.
We write $\Mod A$ for the category of all $A$-modules and $\md A$ for the full subcategory of
finite-dimensional $A$-modules. We write $\dm$ for dimension over $k$, $\rad E$ for the Jacobson radical
of an algebra $E$, and $X^{(I)}$ for the direct sum of copies of $X$ indexed by
a set $I$. Tensor products without a subscript are over $k$; for an $A$-module $P$ and a vector
space $V$, the left $A$-action on $P\otimes V$ is given by $a(p\otimes v)=ap\otimes v$.

For arrows $\alpha:i\to j$ and $\beta:j\to l$ in a quiver,
the product $\beta\alpha$ traverses $\alpha$ first, so that left modules over
a path algebra identify with covariant quiver representations.
We write $\Hom_{K_2}$ and $\End_{K_2}$ for Hom and End
over $kK_2$. A module over a $k$-linear category is a covariant $k$-linear
functor to vector spaces. For a representation $X$ of a quiver, or a
module $X$ over a $k$-linear category, we put $\supp X=\{i:X(i)\ne0\}$.
\end{conventions}

\begin{formalization}
\def\UrlBreaks{\do\/\do\-\do\.}
Theorem~A and the full statement of Corollary~\ref{cor:dimension}
have been formalized in Lean~4 using Mathlib. The source is available at
\url{https://github.com/haruhisa-enomoto/third-brauer-thrall-conjecture}.
Rather than assuming the cited results of Bongartz and G\"obel--Simson as
additional axioms, the development proves the existence of the quotients,
pairs of modules on coverings, and walks needed in
Theorem~\ref{thm:finite-input} (using alternative arguments at the two
steps noted in Remark~\ref{rem:formalization-arguments}), establishes the
numerical bound of Proposition~\ref{prop:global-refinement} directly
without the bounded-bimodule assertion of
Lemma~\ref{lem:quantitative-input}, and proves
Theorem~\ref{thm:rigid-family} for arbitrary infinite fields, which
suffices for our application over algebraically closed fields. In the walk
construction it uses a distinguished sink in place of a source;
Lemma~\ref{lem:L2-identification} and
Appendix~\ref{sec:simson-counterexample} are not formalized.
\end{formalization}

\begin{aiuse}
The research was done under the author's direction using GPT-6-Astra,
which found and developed the proof presented here and wrote the earlier
drafts of the manuscript. The present text was revised with Claude
Opus~5.5. The Lean formalization was carried out mainly by Claude
Opus~5.5. GPT-6-Astra provided mathematical advice for steps that Claude
Opus~5.5 could not complete and performed the final audit.
The author directed the research and revisions and is
responsible for the mathematical content and the final manuscript.
\end{aiuse}

\section{Kronecker modules and nilpotent ideals}
\label{sec:extension}
Throughout this section, $k$ is any field and $A$ is any finite-dimensional
$k$-algebra. We describe the simple summands of arbitrary Kronecker
modules and record how retractions on Hom spaces with nilpotent kernels
preserve indecomposability and reflect isomorphisms.

\subsection{Kronecker modules without simple injective summands}
We first describe $\D$ and $\C$ in terms of arrow maps.
A Kronecker module is written $T=(V,W;a,b)$, with $a,b:V\to W$;
a morphism $(g,h):T\to(V',W';a',b')$ satisfies $ha=a'g$ and $hb=b'g$.
For any Kronecker module $T=(V,W;a,b)$, the submodule
$(0,W;0,0)\cong S_-\otimes W$ and the quotient $(V,0;0,0)\cong S_+\otimes V$
fit into the canonical short exact sequence
\begin{equation}\label{eq:canonical-sequence}
 0\longrightarrow S_-\otimes W\longrightarrow T
 \longrightarrow S_+\otimes V\longrightarrow0.
\end{equation}
Moreover, $\Hom_{K_2}(S_+,T)$ identifies with $\ker a\cap\ker b$:
choosing a complement $V_1$ to $\ker a\cap\ker b$ in $V$ gives a
decomposition $T\cong(V_1,W;a|_{V_1},b|_{V_1})\oplus(\ker a\cap\ker b,0;0,0)$,
where the second summand is a direct sum of copies of $S_+$. Hence
\begin{equation}\label{eq:domain}
 T\in\D\quad\Longleftrightarrow\quad\ker a\cap\ker b=0.
\end{equation}

The following definition describes the subcategory $\C$ announced in
Section~\ref{sec:outline} in terms of the arrow maps.
\begin{definition}\label{def:C}
Let $\C$ be the full subcategory of $\D$ consisting of the modules
$T=(V,W;a,b)$ with $a(V)+b(V)=W$. Write
$\C_{\mathrm{fin}}=\C\cap\md kK_2$.
\end{definition}
For $T\in\D$, choosing a complement to $a(V)+b(V)$ in $W$ gives a
decomposition
\begin{equation}\label{eq:simple-splitting}
 T\cong T_0\oplus S_-^{(I)},\qquad T_0\in\C.
\end{equation}
Here $T_0=(V,a(V)+b(V);a,b)$ and $I$ indexes a basis of the chosen complement.
Conversely, if $T=T'\oplus S_-$ with $T'=(V',W';a',b')$, then
$W=W'\oplus k$ and $a(V)+b(V)=a'(V')+b'(V')\subseteq W'$, so
$a(V)+b(V)\ne W$. Thus $\C$ consists precisely of the Kronecker modules
having neither $S_+$ nor $S_-$ as a direct summand. In particular, every infinite-dimensional indecomposable Kronecker
module belongs to $\C$.

A morphism in $\D$ need not be determined by its component at vertex~$1$ alone, since
every endomorphism of $S_-$ vanishes at vertex~$1$. For modules in $\C$, however,
evaluation at either vertex is injective. We record this fact together
with the closure properties and finite-dimensional submodules of these
categories.

\begin{lemma}\label{lem:Kronecker}
The subcategories $\D$ and $\C$ have the following properties.
\begin{enumerate}
\item Both vertex evaluations $(g,h)\mapsto g$ and $(g,h)\mapsto h$ are
injective on $\Hom_{K_2}(T,U)$ for $T,U\in\C$.
\item $\D$ is closed under submodules, direct sums and filtered colimits,
and $\C$ is closed under finite direct sums.
\item Every object of $\D$ is the union of its finite-dimensional
submodules. For $T=(V,W;a,b)\in\C$, the submodules
\[
 T_{V_0}=(V_0,a(V_0)+b(V_0);a,b),\qquad V_0\subseteq V,\quad\dm V_0<\infty,
\]
are cofinal among its finite-dimensional submodules and belong to $\C$.
\end{enumerate}
\end{lemma}
\begin{proof}
(1) Let $(g,h):T\to U$ be a morphism between $T=(V,W;a,b)$ and $U=(V',W';a',b')$ in $\C$. If $g=0$, then $ha=a'g=0$ and $hb=b'g=0$, so $h=0$ because $a(V)+b(V)=W$.
Conversely, if $h=0$, then $a'g=ha=0$ and $b'g=hb=0$, so $\im g\subseteq\ker a'\cap\ker b'=0$ by \eqref{eq:domain} and hence $g=0$.

(2) A Kronecker module $T=(V,W;a,b)$ satisfies \eqref{eq:domain} if and only if
the map $\binom ab:V\to W\oplus W$ is injective. This condition is inherited by submodules and
preserved by direct sums and filtered colimits, since filtered colimits
of vector spaces are exact. The additional condition $a(V)+b(V)=W$ in Definition~\ref{def:C} is equivalent to the
surjectivity of $(a\ b):V\oplus V\to W$, which is preserved by finite
direct sums.

(3) Since the algebra $kK_2$ is finite-dimensional, every finitely generated
submodule is finite-dimensional, so every object of $\D$ is the union of its finite-dimensional submodules. For $T=(V,W;a,b)\in\C$, each submodule $T_{V_0}$ belongs to $\C$ by definition. Given any finite-dimensional submodule $U=(V_1,W_1;a,b)\subseteq T$, any finite spanning set of $W_1$ lies in $a(V_0)+b(V_0)$ for a sufficiently large
finite-dimensional subspace $V_0\subseteq V$; enlarging $V_0$ so that $V_1\subseteq V_0$ gives $U\subseteq T_{V_0}$, which proves cofinality.
\end{proof}

\subsection{Nilpotent kernels and indecomposability}
\label{sec:nilpotent-kernels}
To show that a functor $\Psi$ preserves indecomposability and reflects
isomorphism classes, we will recover morphisms between Kronecker modules
from morphisms between their images under $\Psi$.
Ardizzoni and Menini \cite[Definition~1.1]{ArdizzoniMenini} call a functor
\emph{heavily separable} if it admits natural retractions on Hom spaces
preserving composition. We use $k$-linear retractions of this kind whose
endomorphism kernels satisfy a uniform nilpotence bound.

\begin{definition}\label{def:retraction}
Let $\Psi$ be a $k$-linear functor from a full subcategory of $\Mod kK_2$
to $\Mod A$, and let $\mathcal X$ be a full subcategory
of its domain. We call $k$-linear maps
\[
 \rho_{X,Y}:\Hom_A(\Psi(X),\Psi(Y))\longrightarrow\Hom_{K_2}(X,Y),\qquad
 X,Y\in\mathcal X,
\]
\emph{retractions for $\Psi$ on $\mathcal X$} if they preserve composition
and satisfy $\rho_{X,Y}(\Psi(u))=u$ for every morphism $u:X\to Y$. Thus
each $\rho_{X,Y}$ is a retraction of the linear map
$\Hom_{K_2}(X,Y)\to\Hom_A(\Psi(X),\Psi(Y))$ induced by $\Psi$.
For an integer $n\ge1$, these retractions have \emph{nilpotence bound $n$} if
\begin{equation}\label{eq:nilpotence-bound}
 (\ker\rho_{X,X})^n=0\qquad\text{for every }X\in\mathcal X.
\end{equation}
\end{definition}
Taking $u=1_X$ shows that $\rho_{X,X}(1_{\Psi(X)})=1_X$, so each $\rho_{X,X}$ is a unital $k$-algebra homomorphism and $\ker\rho_{X,X}$ is a two-sided ideal of $\End_A(\Psi(X))$. Moreover, composition
with maps in the image of $\Psi$ shows that $\rho_{X,Y}$ is natural in both variables.

An infinite-dimensional indecomposable module need not have a local
endomorphism algebra: for instance, $(k[t],k[t];1,t)$ has endomorphism algebra $k[t]$.
Instead, we detect indecomposability through the absence of nontrivial idempotents.
\begin{lemma}\label{lem:retraction-consequences}
Let $\Psi:\C\to\Mod A$ be a $k$-linear functor admitting retractions
with nilpotence bound $n\ge1$. Then $\Psi$ preserves indecomposability
and reflects isomorphism classes.
\end{lemma}
\begin{proof}
Let $T\in\C$ be indecomposable. Since $\rho_{T,T}(\Psi(1_T))=1_T\ne0$, we have $\Psi(T)\ne0$. Let $e\in\End_A(\Psi(T))$ be an
idempotent. Because $\rho_{T,T}$ is a unital algebra homomorphism and $T$ is indecomposable, the idempotent
$\rho_{T,T}(e)\in\End_{K_2}(T)$ is either $0$ or $1_T$. Consequently, either $e$ or $1_{\Psi(T)}-e$ lies in $\ker\rho_{T,T}$ and is therefore a
nilpotent idempotent and hence zero. Thus $\Psi(T)$ is indecomposable.
Finally, if $T,U\in\C$ and $f:\Psi(T)\to\Psi(U)$ is an isomorphism, then
$\rho_{T,U}(f):T\to U$ and $\rho_{U,T}(f^{-1}):U\to T$ are mutually inverse morphisms in $\C$, so $T\cong U$.
\end{proof}

\section{Bongartz's representation embeddings}
\label{sec:constructions}
In Section~\ref{sec:reduction}, $k$ is any field and $A$ is a
representation-infinite finite-dimensional $k$-algebra. In the remaining
subsections, $k$ is algebraically closed and $A$ denotes a connected basic
minimal representation-infinite finite-dimensional $k$-algebra.
We define the three functors used in Bongartz's representation embedding
theorem \cite[Theorem~6]{Bongartz} on modules of arbitrary dimension.
Theorem~\ref{thm:finite-input} supplies the required algebras, modules and walks, and
Section~\ref{sec:control} proves the nilpotence needed to preserve
indecomposability.

\subsection{Reduction to minimal representation-infinite algebras}
\label{sec:reduction}
Recall that an algebra is called \emph{minimal representation-infinite}
if it is representation-infinite and all of its proper quotients are
representation-finite.
The following lemma records the functor used in the reduction of
\cite[\S5.2]{Bongartz} to transfer
embeddings back to $A$ in arbitrary dimension.

\begin{lemma}\label{lem:reduction}
Let $A$ be a representation-infinite finite-dimensional algebra over a
field $k$. Choose a full idempotent $e$ (so $AeA=A$) such that $A_0=eAe$ is basic,
and a minimal representation-infinite quotient $B$ of $A_0$.
Then $B$ is basic and connected, $\dm B\le\dm A$, and there is an
exact fully faithful functor $\Phi:\Mod B\to\Mod A$ preserving
finite-dimensional modules, arbitrary direct sums, and filtered
colimits. For every infinite-dimensional $B$-module $X$,
$\dm\Phi(X)=\dm X$.
\end{lemma}
\begin{proof}
Since $A$ and $A_0=eAe$ are Morita equivalent, $A_0$ is representation-infinite, and any quotient $B$ of $A_0$ of least dimension among its
representation-infinite quotients is minimal
representation-infinite with $\dm B\le\dm A_0\le\dm A$. Any minimal
representation-infinite quotient $B$ of $A_0$ is basic because $A_0$ is basic, and
connected because otherwise one of its connected factors would be a proper
representation-infinite quotient. Let $\Phi:\Mod B\to\Mod A$ be inflation from $B$ to $A_0$
followed by the Morita equivalence $Ae\otimes_{A_0}-$.
Both inflation and $Ae\otimes_{A_0}-$ are exact and fully faithful, preserve finite-dimensional modules, and
commute with arbitrary direct sums and filtered colimits. Since the inverse Morita functor is $e(-)$, every $B$-module $X$ satisfies
\[
 \dm X=\dm e\Phi(X)\le\dm\Phi(X)\le(\dm Ae)\dm X,
\]
and hence $\dm\Phi(X)=\dm X$ whenever $\dm X$ is infinite.
\end{proof}

In the remainder of this section, we write $A$ for the reduced algebra $B$
of Lemma~\ref{lem:reduction}, as stipulated at the start of this section.

\subsection{Quotients of projective modules by submodules of the socle}
\label{sec:quotients}
In the first construction, a Kronecker module is sent to a quotient of a direct sum of copies of one indecomposable projective module over an algebra $\Lambda$, with relations defined by the two arrow maps.

\begin{setting}\label{set:quotient}
Throughout this subsection, let $\Lambda$ be a basic finite-dimensional
algebra over the algebraically closed field $k$, and put $J=\rad\Lambda$.
Fix a complete set $\{e_z\}$ of primitive orthogonal idempotents of $\Lambda$.
Choose two idempotents $e_x,e_y$ from this set (which may coincide) and
linearly independent elements $s,t\in e_xJe_y$ satisfying $Js=Jt=sJ=tJ=0$.
\end{setting}

With the socle elements $s,t$ of Setting~\ref{set:quotient}, the functor
of \cite[\S3.1]{Bongartz} takes the following form, as used in the proof
of \cite[Corollary~2]{Bongartz}.
\begin{construction}[{\cite[\S3.1]{Bongartz}}]
\label{con:quotient}
For $T=(V,W;a,b)\in\D$, put
\begin{equation}\label{eq:quotient-construction}
 F_1(T)=(\Lambda e_y\otimes W)/
 \{s\otimes a(v)+t\otimes b(v):v\in V\}.
\end{equation}
On a morphism $(g,h)$, let $F_1(g,h)$ be the map induced by $1\otimes h$.
\end{construction}
To see that the relation space $\{s\otimes a(v)+t\otimes b(v):v\in V\}$ is a $\Lambda$-submodule of $\Lambda e_y\otimes W$, note that it
is annihilated by $J$ and satisfies $e_zs=\delta_{zx}s$ and $e_zt=\delta_{zx}t$ for
each primitive idempotent $e_z$; since $\Lambda/J\cong\prod_z k$, it is
stable under $\Lambda$ and lies in
the socle of $\Lambda e_y\otimes W$. For any morphism $(g,h):T\to(V',W';a',b')$, the relations
$ha=a'g$ and $hb=b'g$ imply that $1\otimes h$ maps the relation submodule for $T$ into that for $(V',W';a',b')$, so $F_1(g,h)$ is well defined. Moreover, since $s,t\in J$, the relation submodule is contained in
$J(\Lambda e_y\otimes W)$, so the top of $F_1(T)$ is identified with $S_y\otimes W$, where
$S_y=\Lambda e_y/Je_y$; we use this identification of the induced map on tops in
Proposition~\ref{prop:quotient-control}.

\begin{lemma}\label{lem:F1-properties}
The functor $F_1:\D\to\Mod\Lambda$ is exact and faithful, commutes
with arbitrary direct sums and filtered colimits, and preserves
finite-dimensional modules. For infinite-dimensional $T\in\D$,
$\dm F_1(T)=\dm T$.
\end{lemma}
\begin{proof}
For $T=(V,W;a,b)\in\D$, define $j_T:V\to\Lambda e_y\otimes W$ by
$j_T(v)=s\otimes a(v)+t\otimes b(v)$. Since $s$ and $t$ are linearly
independent and $T\in\D$, we have $\ker j_T=\ker a\cap\ker b=0$ by \eqref{eq:domain}, so $j_T$ is injective.
Viewing $V\cong S_x\otimes V$ as a semisimple $\Lambda$-module via the simple
module $S_x=\Lambda e_x/Je_x$, we may regard $j_T$ as a $\Lambda$-module homomorphism with cokernel $F_1(T)$ (so the formula \eqref{eq:quotient-construction} defines a right exact functor on $\Mod kK_2$).
For a short exact sequence
$0\to T'\to T\to T''\to0$ in $\D$, the commutative diagram
\[
\begin{tikzcd}[column sep=small]
0\arrow[r]&V'\arrow[r]\arrow[d,"j_{T'}"]&V\arrow[r]\arrow[d,"j_T"]&V''\arrow[r]\arrow[d,"j_{T''}"]&0\\
0\arrow[r]&\Lambda e_y\otimes W'\arrow[r]&\Lambda e_y\otimes W\arrow[r]&\Lambda e_y\otimes W''\arrow[r]&0
\end{tikzcd}
\]
has exact rows and injective vertical maps, so the induced sequence of cokernels
$0\to F_1(T')\to F_1(T)\to F_1(T'')\to0$ is exact.
Thus $F_1$ is exact on $\D$. If $(g,h):T\to U$ is a morphism in $\D$ with $F_1(g,h)=0$, then its induced map $1_{S_y}\otimes h$ on tops vanishes, so
$h=0$ and hence $g=0$ by \eqref{eq:domain} for $U$;
thus $F_1$ is faithful. Since tensor products and cokernels commute with
arbitrary direct sums and filtered colimits, so does $F_1$. Finally, $F_1$ preserves
finite-dimensional modules, and the inequalities
\[
 \dm W\le\dm F_1(T)\le(\dm\Lambda e_y)\dm W
\]
together with the embedding $\binom ab:V\hookrightarrow W\oplus W$ show that $\dm F_1(T)=\dm W=\dm T$ whenever $T\in\D$ is infinite-dimensional.
\end{proof}

\subsection{Galois coverings and push-down functors}
\label{sec:coverings}
Fix a complete set $\{e_q\}$ of primitive orthogonal idempotents of the basic
algebra $A$, indexed by a set of vertices $q$, and let $\mathcal A$ be the $k$-linear
category with these vertices as objects and morphism spaces
$\mathcal A(q,q')=e_{q'}Ae_q$; we identify $A$-modules with
$\mathcal A$-modules. In the second construction, a larger category $R$
carries an action of a group $G$ that is free on objects, and a functor
$\pi:R\to\mathcal A$ identifies each $G$-orbit of objects with a single object
of $\mathcal A$. Kronecker modules are first embedded fully
faithfully into modules over $R$ supported on finitely many objects, and the
spaces above each object of $\mathcal A$ are then summed.

\begin{definition}[{\cite[3.1]{GabrielUniversal}}]\label{def:covering}
A small $k$-linear category $R$ is \emph{locally bounded} if distinct
objects are nonisomorphic, each endomorphism algebra is local, all
morphism spaces are finite-dimensional, and for every object $x$ the set
$\{y:R(x,y)\ne0\text{ or }R(y,x)\ne0\}$ is finite.
A \emph{Galois covering with group $G$} is a $k$-linear functor
$\pi:R\to\mathcal A$ from a locally bounded category $R$, surjective on
objects, together with an action
of $G$ by $k$-linear automorphisms, such that the action on objects is
free, $\pi g=\pi$, each fibre is a $G$-orbit, and both maps
\begin{equation}\label{eq:covering-lifts}
 \bigoplus_{\pi y=q}R(x,y)\longrightarrow\mathcal A(\pi x,q),
 \qquad
 \bigoplus_{\pi y=q}R(y,x)\longrightarrow\mathcal A(q,\pi x)
\end{equation}
induced by $\pi$ are isomorphisms for every $x,q$.
\end{definition}

For the rest of this subsection, fix a Galois covering
$\pi:R\to\mathcal A$ with group $G$ as in Definition~\ref{def:covering}.
We allow arbitrary vertex spaces, but require the support to be finite.
Its size will bound the nilpotence exponent in Section~\ref{sec:control}.
\begin{construction}\label{con:pushdown}
For an $R$-module $X$ of finite support, set
\[
 (\pi_*X)(q)=\bigoplus_{\pi x=q}X(x).
\]
Given a morphism $u:q\to q'$ in $\mathcal A$ and an object $x$ with
$\pi x=q$, the first isomorphism in \eqref{eq:covering-lifts} allows us to write
$u=\sum_{\pi y=q'}\pi(u_{y,x})$ with unique $u_{y,x}\in R(x,y)$. We define
$(\pi_*X)(u):(\pi_*X)(q)\to(\pi_*X)(q')$ to be the linear map whose $(y,x)$-block is $X(u_{y,x}):X(x)\to X(y)$, and on morphisms of $R$-modules we define $\pi_*$ by taking the direct sum of the
vertex maps. This defines the \emph{push-down functor} $\pi_*$.
\end{construction}
By the uniqueness in \eqref{eq:covering-lifts}, these block matrices
respect composition.

We consider modules supported on a finite full subcategory $\mathcal S\subset R$
satisfying
\begin{equation}\label{eq:extension-by-zero}
 R(x,j)\circ R(i,x)=0\qquad(i,j\in\mathcal S,\ x\notin\mathcal S).
\end{equation}
Write $\Gamma=\bigoplus_{i,j\in\mathcal S}R(i,j)$, with multiplication induced
by composition in $R$ and zero on noncomposable factors, so that $\Gamma$-modules identify with modules over $\mathcal S$. Condition
\eqref{eq:extension-by-zero} ensures that extending a $\Gamma$-module by
zero outside $\mathcal S$ yields an $R$-module, and that extension by zero is an exact fully faithful
functor from $\Gamma$-modules to $R$-modules.

We use the construction by extensions of Ringel
\cite[1.5, Lemma]{RingelSpecies}. Its form for modules of arbitrary dimension is
also used by G\"obel and Simson \cite[Proposition~3.1]{GoebelSimsonPrinjective}
and Simson \cite[proof of Theorem~1.2]{SimsonRealization}.
\begin{construction}\label{con:cover}
Choose finite-dimensional $\Gamma$-modules $M_-,M_+$ with disjoint
supports and $\End_\Gamma(M_-)=\End_\Gamma(M_+)=k$, and two linearly independent
classes in $\Ext^1_\Gamma(M_+,M_-)$. Fix a bound quiver presentation
of $\Gamma$. Represent each class by an extension with vertex spaces
$M_-(i)\oplus M_+(i)$, and write $z_j(\alpha):M_+(i)\to M_-(l)$
for the off-diagonal arrow entries of the $j$-th extension, for
$\alpha:i\to l$ and $j=1,2$.
For any Kronecker module $T=(V,W;a,b)$, define $L_2(T)$ by
\[
 L_2(T)(i)=(M_-(i)\otimes W)\oplus(M_+(i)\otimes V),\qquad
 L_2(T)(\alpha)=
 \begin{pmatrix}
 M_-(\alpha)\otimes1_W&z_1(\alpha)\otimes a+z_2(\alpha)\otimes b\\
 0&M_+(\alpha)\otimes1_V
 \end{pmatrix}.
\]
For a morphism $(g,h):T\to T'$, we define $L_2(g,h)$ at each vertex $i\in\mathcal S$ by the diagonal blocks $1_{M_-(i)}\otimes h$ and $1_{M_+(i)}\otimes g$. Viewing $L_2(T)$ as an $R$-module via extension by zero outside $\mathcal S$, we put
$F_2=\pi_*L_2:\Mod kK_2\to\Mod A$.
\end{construction}
To check that $L_2(T)$ satisfies the defining relations of $\Gamma$, observe that in any product of the upper triangular matrices $L_2(T)(\alpha)$, each term
in the off-diagonal block contains a single off-diagonal factor.
For any relation of $\Gamma$, the off-diagonal block of the evaluated relation is therefore the sum of the corresponding
blocks for the two chosen extensions, tensored with $a$ and $b$ respectively.
Both vanish because those extensions are $\Gamma$-modules, and the diagonal
blocks satisfy the relations because $M_-$ and $M_+$ do, so $L_2(T)$
is a $\Gamma$-module in arbitrary dimension.

\begin{lemma}\label{lem:L2}
The functor $L_2$ is exact and fully faithful. The functor $F_2$ is
exact and faithful. Both commute with direct sums and filtered colimits
and preserve finite-dimensional modules. For infinite-dimensional $T$,
$\dm F_2(T)=\dm T$.
\end{lemma}
\begin{proof}
Let $T=(V,W;a,b)$ and $T'=(V',W';a',b')$ be Kronecker modules. Since $M_-$ and $M_+$ have disjoint supports, at each vertex $i\in\mathcal S$ at most one of $M_-(i)$ and $M_+(i)$ is nonzero. Any $\Gamma$-module homomorphism $L_2(T)\to L_2(T')$ therefore decomposes into maps $M_-\otimes W\to M_-\otimes W'$ and $M_+\otimes V\to M_+\otimes V'$, and the diagonal blocks of the arrow equations show that these two components are $\Gamma$-module homomorphisms.
For a finite-dimensional $\Gamma$-module $M$ with
$\End_\Gamma(M)=k$, every homomorphism $M\otimes U\to M\otimes U'$ has the form
$1_M\otimes h$ for a unique linear map $h:U\to U'$: indeed,
$M\otimes U$ is a direct sum of copies of $M$, and since $M$ is finitely generated, each copy maps into a finite direct sum of copies of $M$ via scalar endomorphisms of $M$, which define the map $h$. Thus the two components of our morphism are $1_{M_-}\otimes h$
and $1_{M_+}\otimes g$ for unique linear maps $h:W\to W'$ and $g:V\to V'$.
The remaining arrow equations are
\[
 z_1(\alpha)\otimes(ha-a'g)+
 z_2(\alpha)\otimes(hb-b'g)=0
 \quad\text{for every arrow }\alpha.
\]
Since the two families $z_1,z_2$ represent linearly independent extension classes, they are linearly independent in the
finite-dimensional space
\[
 \bigoplus_{\alpha:i\to l}\Hom_k(M_+(i),M_-(l)).
\]
Applying linear functionals that separate $z_1$ and $z_2$ yields
$ha=a'g$ and $hb=b'g$, which proves that $L_2$ is fully faithful.
Since $L_2$ is defined at each vertex by tensoring with finite-dimensional spaces and taking direct sums, it is exact, preserves finite-dimensional modules, and commutes with arbitrary direct sums and filtered colimits.
Extending by zero to $R$ and applying $\pi_*$ preserve these properties, and $\pi_*$ is exact and
faithful because it sums vertex spaces and maps over each fibre. Finally,
\[
 \dm F_2(T)=(\dm M_-)\dm W+(\dm M_+)\dm V,
\]
and since $\dm M_-$ and $\dm M_+$ are positive integers, $\dm F_2(T)=\dm T$ whenever $T$ is infinite-dimensional.
\end{proof}

In the covering case of Theorem~\ref{thm:finite-input}, Bongartz supplies
a functor on finite-dimensional modules. The next lemma identifies that functor with the finite restriction
of $L_2$, so that Construction~\ref{con:cover} extends the specified
functor to arbitrary dimension.
\begin{lemma}\label{lem:L2-identification}
Let $E:\md kK_2\to\md\Gamma$ be an exact fully faithful $k$-linear
functor such that $E(S_-)$ and $E(S_+)$ have disjoint supports.
Put $M_-=E(S_-)$ and $M_+=E(S_+)$. Represent the images under $E$ of
the extension classes with middle terms $(k,k;1,0)$ and $(k,k;0,1)$
by extensions, and let $z_1,z_2$ be their off-diagonal arrow entries
as in Construction~\ref{con:cover}.
Then these data satisfy the hypotheses of that construction, and $E$
is naturally isomorphic to the restriction to $\md kK_2$ of the
resulting functor $L_2$.
\end{lemma}
\begin{proof}
Since $E$ is fully faithful, we have $\End_\Gamma(M_-)=\End_\Gamma(M_+)=k$.
The $k$-linear map $\Ext^1_{K_2}(S_+,S_-)\to\Ext^1_\Gamma(M_+,M_-)$ induced by the exact $k$-linear functor
$E$ is injective, since if the image of an extension splits in $\md\Gamma$, full
faithfulness lifts its splitting to $\md kK_2$. Thus the images of the two basis extensions with middle terms $(k,k;1,0)$ and $(k,k;0,1)$ are linearly independent.
For any finite-dimensional Kronecker module $T=(V,W;a,b)$, applying $E$ to \eqref{eq:canonical-sequence} yields an exact sequence
$0\to M_-\otimes W\to E(T)\to M_+\otimes V\to0$.
Because $M_-$ and $M_+$ have disjoint supports, at each vertex $i\in\mathcal S$ at most one of $M_-(i)$ and $M_+(i)$ is nonzero. The exact sequence therefore gives a canonical identification
$E(T)(i)\cong(M_-(i)\otimes W)\oplus(M_+(i)\otimes V)=L_2(T)(i)$ at every vertex $i\in\mathcal S$, and every vertex linear map
$M_+(i)\otimes V\to M_-(i)\otimes W$ vanishes. Consequently, there are no nonzero coboundaries on the vertex spaces, so an extension class in $\Ext^1_\Gamma(M_+\otimes V,M_-\otimes W)$ determines unique off-diagonal arrow maps $M_+(i)\otimes V\to M_-(l)\otimes W$.
Since the extension class of \eqref{eq:canonical-sequence} has coefficients $a$ and $b$
with respect to the two basis extensions, under this identification of vertex spaces the off-diagonal arrow maps of $E(T)$ are
$z_1(\alpha)\otimes a+z_2(\alpha)\otimes b$. Since these identifications commute with the maps induced by any
Kronecker morphism $(g,h)$, they define a natural isomorphism $E\cong L_2|_{\md kK_2}$.
\end{proof}

\subsection{Push-down along a closed walk}
\label{sec:walks}
The third construction likewise sums vertex spaces over fibres: a closed walk in
a quiver $Q$ determines a cycle quiver together with a map from that cycle to $Q$, and
we push representations down along this map.

\begin{definition}\label{def:walk}
A walk in a quiver $Q$ is a composable sequence of arrows and formal inverse arrows,
written in traversal order; it is \emph{closed} if it ends at its
starting vertex, and \emph{reduced} if adjacent letters are
never mutual inverses. A nonempty closed walk is \emph{cyclically reduced}
if this also holds at the join of its last and first letters, and
\emph{primitive} if it is not a proper power of a shorter closed walk.
For a closed walk
$w=c_1^{\varepsilon_1}\cdots c_r^{\varepsilon_r}$ passing through vertices
$q_0,q_1,\ldots,q_r=q_0$, we form a cycle quiver $D_w$ with vertices
$d_0,\ldots,d_{r-1}$ (indices read modulo $r$, often written simply as $i$ for $d_i$), where the $i$-th arrow goes from $d_{i-1}$ to $d_i$ if $\varepsilon_i=1$
and in the reverse direction if $\varepsilon_i=-1$. The quiver map $p:D_w\to Q$ sends
$d_i$ to $q_i$ and the $i$-th arrow to $c_i$; we call $p(e)$ the label
of an arrow $e$ of $D_w$.
\end{definition}

When $w$ is cyclically reduced,
\begin{equation}\label{eq:local-injectivity}
 p\text{ is injective on the arrows of $D_w$ with a fixed source, and on
 those with a fixed target,}
\end{equation}
since two outgoing arrows with the same label, or two incoming arrows
with the same label, would correspond to adjacent mutually inverse letters in $w$
(possibly at the join).

\begin{setting}\label{set:walk}
For the rest of this subsection, fix a presentation $A=kQ/I$ as a bound
quiver algebra and a primitive cyclically reduced closed walk $w$ whose
cycle $D_w$ has a source $u$ and $r\ge2$ vertices. Write $p:D_w\to Q$
for the map of Definition~\ref{def:walk}, and fix an ordering of the two
arrows out of $u$.
\end{setting}

\begin{construction}[{\cite[proof of Theorem~7]{Bongartz}}]\label{con:walk}
For $T=(V,W;a,b)\in\Mod kK_2$, let $L_3(T)$ be the representation of $D_w$ with $V$ at $u$, $W$ at all other
vertices of $D_w$, the maps $a,b$ on the ordered arrows out of $u$, and the identity
map $1_W$ on the remaining arrows; a morphism $(g,h)$ acts as $g$ at $u$ and as $h$ at all other vertices.
For a representation $X$ of $D_w$ and a vertex $q$ of $Q$, set
\[
 (p_*X)(q)=\bigoplus_{p(i)=q}X(i).
\]
For an arrow $c:q\to q'$ of $Q$, we define $(p_*X)(c)$ so that its block from $X(i)$ to $X(j)$
is $X(e)$ if there is an arrow $e:i\to j$ in $D_w$ with label $c$, and zero
otherwise; such an arrow is unique by \eqref{eq:local-injectivity}, so
each block is well defined. On morphisms, $p_*$ takes the direct sum of the
vertex maps. The composite $F_3=p_*L_3$ is defined on
$\Mod kK_2$ and takes values in representations of $Q$.
\end{construction}
Theorem~\ref{thm:finite-input}(3) supplies a walk for which the images
of all finite-dimensional Kronecker modules are annihilated by $I$.
The next lemma extends this property to arbitrary modules.
\begin{lemma}\label{lem:F3-properties}
In Setting~\ref{set:walk}, the functor $L_3$ is exact and fully faithful,
and $L_3(S_+)$ and $L_3(S_-)$ are nonzero with disjoint supports.
Suppose that $I$ annihilates $F_3(T)$ for every finite-dimensional
Kronecker module $T$. Then $F_3:\Mod kK_2\to\Mod A$ is exact and
faithful, commutes with arbitrary direct sums and filtered colimits,
and preserves finite-dimensional modules. For infinite-dimensional $T$,
$\dm F_3(T)=\dm T$.
\end{lemma}
\begin{proof}
Because the identity arrows of $L_3(T)$ connect all vertices of $D_w$ other than $u$, any morphism $L_3(T)\to L_3(V',W';a',b')$
must have the same component $h$ at all those vertices, and the remaining two commutativity conditions at the arrows out of $u$ are precisely $ha=a'g$ and $hb=b'g$. Thus $L_3$ is fully
faithful; it is also exact because each vertex space is either $V$ or $W$. Moreover,
$L_3(S_+)$ and $L_3(S_-)$ are nonzero and have disjoint supports $\{u\}$
and $(D_w)_0\setminus\{u\}$, just as for $L_2$ in Construction~\ref{con:cover}.

Since the push-down $p_*$ takes finite direct sums of vertex
spaces and maps, it is exact and faithful, and both $L_3$ and $p_*$ commute
with arbitrary direct sums and filtered colimits and preserve finite-dimensional modules.
Because every Kronecker module is the union of its finite-dimensional
submodules and $F_3=p_*L_3$ commutes with filtered colimits, the ideal
$I$ annihilates $F_3(T)$ for every $T\in\Mod kK_2$, so $F_3$ takes values in $\Mod A$ and is exact and faithful.
Finally, $\dm F_3(T)=\dm V+(r-1)\dm W=\dm T$ whenever $T$ is
infinite-dimensional.
\end{proof}

\subsection{Existence of the embeddings}
\label{sec:finite-input}
Over the algebraically closed field $k$, Bongartz \cite[Theorem~6 and its proof]{Bongartz} shows that at least one
of the three preceding constructions yields a representation embedding on
$\D_{\mathrm{fin}}$. While \cite[Theorem~6]{Bongartz} states the existence of an embedding, the free covering group,
finite full subcategory and primitive closed walk in the alternatives
below come from the proof of \cite[Theorem~7]{Bongartz} in \cite[\S5.4]{Bongartz}.
The subcategory $\mathcal S$ chosen there has the property that a nonzero
composition between its objects cannot pass through an object outside
$\mathcal S$; this gives \eqref{eq:extension-by-zero}.
The disjoint supports in alternative (2) follow from
\cite[Lemma~5]{Bongartz}, as we explain after the theorem.

\begin{theorem}[{\cite[Theorem~6 and its proof]{Bongartz}}]
\label{thm:finite-input}
Let $A$ be a connected basic minimal representation-infinite algebra
over an algebraically closed field. Then at least one of the following
holds.
\begin{enumerate}
\item $A$ has elements $s,t$ as in Setting~\ref{set:quotient}, with
$\Lambda=A$, such that $F_1$ from Construction~\ref{con:quotient} is a
representation embedding $\D_{\mathrm{fin}}\to\md A$
\cite[Theorem~1 and Corollary~2]{Bongartz}.
\item There are a Galois covering $\pi:R\to\mathcal A$ whose group
$G$ is a free group, a finite full subcategory $\mathcal S$ of $R$
satisfying \eqref{eq:extension-by-zero}, and
an exact fully faithful $k$-linear functor $E:\md kK_2\to\md\Gamma$
such that $E(S_+)$ and $E(S_-)$ are nonzero with disjoint supports and
$\pi_*E|_{\D_{\mathrm{fin}}}$, with $E$ extended by zero from
$\mathcal S$ to $R$, is a representation embedding
$\D_{\mathrm{fin}}\to\md A$
\cite[Theorem~3, Lemma~5, and proof of Theorem~7]{Bongartz}.
\item There are a presentation $A=kQ/I$ and a walk as in
Setting~\ref{set:walk} such that $F_3$ from Construction~\ref{con:walk}
sends finite-dimensional Kronecker modules to $\md A$ and
$F_3|_{\D_{\mathrm{fin}}}$ is a
representation embedding $\D_{\mathrm{fin}}\to\md A$
\cite[proof of Theorem~7]{Bongartz}.
\end{enumerate}
\end{theorem}
In alternative (1), the construction in \cite[Corollary~2]{Bongartz}
initially uses a quotient $\Lambda$ of $A$. Since the modules $(k,k;\lambda,1)$
($\lambda\in k$) are pairwise non-isomorphic indecomposables in
$\D_{\mathrm{fin}}$, the existence of a representation embedding $\D_{\mathrm{fin}}\to\md\Lambda$
implies that $\Lambda$ is representation-infinite; by the minimality of $A$, we have $\Lambda=A$.

For alternative (2), the proof of \cite[Theorem~7]{Bongartz} chooses
$\mathcal S$ so that $\Gamma$ is tame concealed and verifies the
representation-embedding property after push-down.
By \cite[Lemma~5 and its proof]{Bongartz}, there is a source or sink $x$
and, in the sink case, an exact sequence of $\Gamma$-modules
$0\to P_x\to N_x\to U_x\to0$, where $P_x$ is simple projective,
$\dm N_x(x)=1$, and
\[
 \End_\Gamma(P_x)=\End_\Gamma(U_x)=k,\qquad
 \Hom_\Gamma(P_x,U_x)=\Hom_\Gamma(U_x,P_x)=0,\qquad
 \dm\Ext^1_\Gamma(U_x,P_x)=2.
\]
Since $P_x$ is supported at $x$ and $\dm N_x(x)=1$, exactness at $x$ gives $U_x(x)=0$, so $M_-=P_x$ and $M_+=U_x$
have disjoint supports. When $x$ is a source, the lemma is applied to
$\Gamma^{\mathrm{op}}$: $P_x,N_x,U_x$ are then
$\Gamma^{\mathrm{op}}$-modules, and we take the $\Gamma$-modules
$M_-=U_x^*$ and $M_+=P_x^*$, where $(-)^*=\Hom_k(-,k)$.
Duality preserves disjointness and identifies
$\Ext^1_\Gamma(P_x^*,U_x^*)$ with
$\Ext^1_{\Gamma^{\mathrm{op}}}(U_x,P_x)$.
The functor $E$ used in \cite[proof of Theorem~7]{Bongartz} is the
finite-dimensional restriction of the extension functor of Construction~\ref{con:cover} for this pair
$(M_-,M_+)$, so $E(S_-)=M_-$ and $E(S_+)=M_+$.
Any choice of basis of the two-dimensional space $\Ext^1_\Gamma(M_+,M_-)$ gives
the same representation-embedding property after push-down, since a change
of basis replaces the two Kronecker arrow maps by invertible linear
combinations and thus acts by an autoequivalence of $\md kK_2$ preserving $\D_{\mathrm{fin}}$.

For each available alternative $\nu$ of Theorem~\ref{thm:finite-input},
fix its data and write
\begin{equation}\label{eq:chosen-functors}
 F_\nu:\D\longrightarrow\Mod A
\end{equation}
for the functor given by that construction, with $F_2,F_3$ restricted to $\D$.
In alternative (2), Lemma~\ref{lem:L2-identification} supplies the
Construction~\ref{con:cover} data from the specified exact fully faithful
$k$-linear functor $E$, so the finite restriction of $F_2$ is naturally
isomorphic to $\pi_*E|_{\D_{\mathrm{fin}}}$.
Put $F_{\nu,\mathrm{fin}}=F_\nu|_{\D_{\mathrm{fin}}}$. Since every
$T\in\D$ is the union of its finite-dimensional submodules,
\begin{equation}\label{eq:colimit-comparison}
 F_\nu(T)\cong\varinjlim_{\substack{U\subseteq T\\\dm U<\infty}}
 F_{\nu,\mathrm{fin}}(U),
\end{equation}
naturally in $T\in\D$. This identifies $F_\nu$ with the extension
of $F_{\nu,\mathrm{fin}}$ by filtered colimits, uniquely up to natural
isomorphism among extensions preserving filtered colimits. The same holds after
applying $\Phi$ from Lemma~\ref{lem:reduction}; we use it in Step~1
of the proof of Theorem~\ref{thm:embedding}.

\begin{remark}\label{rem:formalization-arguments}
The Lean formalization gives alternative arguments for the pair of modules
with disjoint supports in Theorem~\ref{thm:finite-input}(2) and the closed walk in
Theorem~\ref{thm:finite-input}(3), replacing respectively
\cite[Lemma~5]{Bongartz} and the use of the classification of large sincere
algebras in \cite[proof of Theorem~7]{Bongartz}.
\end{remark}

\section{Homomorphisms between the embedded modules}
\label{sec:control}
Throughout this section, $k$ is algebraically closed and $A$ denotes a
connected basic minimal representation-infinite finite-dimensional
$k$-algebra. Working with Kronecker modules of arbitrary dimension, whose
images under the covering and walk constructions still have finite support
on $R$ and $D_w$, we construct retractions on $\C$ for each of the functors
$F_\nu$ of \eqref{eq:chosen-functors}, with a nilpotence bound determined by
the Loewy length of $A$ or the maximum fibre size of the support.

\subsection{Homomorphisms between quotients by submodules of the socle}
We work in Setting~\ref{set:quotient} and Construction~\ref{con:quotient}.
As noted after Construction~\ref{con:quotient}, the top of
$F_1(T)$ is $S_y\otimes W$, and multiplying lifts of
this top by $s$ and $t$ defines a linear map whose kernel is precisely
$\{(a(v),b(v)):v\in V\}$ (see \eqref{eq:recover-relation} below).
In the proof of \cite[Theorem~1]{Bongartz}, Bongartz uses this observation to recover Kronecker
morphisms from the induced maps on tops; we formulate this recovery as a system of
retractions.

\begin{proposition}\label{prop:quotient-control}
There are retractions
\[
 \rho_{T,U}:\Hom_\Lambda(F_1(T),F_1(U))\longrightarrow\Hom_{K_2}(T,U),
 \qquad T,U\in\D,
\]
for $F_1$ on $\D$.
If $J^n=0$, these retractions have nilpotence bound $n$.
\end{proposition}
\begin{proof}
For any $\Lambda$-module $M$, put $\operatorname{top}_y(M)=e_y(M/JM)$.
Since $sJ=tJ=0$, the map
\[
 d_M:\operatorname{top}_y(M)\oplus\operatorname{top}_y(M)\longrightarrow e_xM,
 \qquad (w_1,w_2)\longmapsto sm_1+tm_2,
\]
where $m_i\in e_yM$ lifts $w_i$, is independent of the choice of lifts.
For $T=(V,W;a,b)\in\D$, under the identification $\operatorname{top}_y(F_1(T))=W$,
the map $d_{F_1(T)}$ is the composite of the injection
$W\oplus W\to\Lambda e_y\otimes W$,
$(w_1,w_2)\mapsto s\otimes w_1+t\otimes w_2$,
with the quotient map onto $F_1(T)$.
Its kernel is therefore
\begin{equation}\label{eq:recover-relation}
 \ker d_{F_1(T)}=\{(a(v),b(v)):v\in V\}.
\end{equation}

Now write $U=(V',W';a',b')$. Given $f:F_1(T)\to F_1(U)$, let
$h:W\to W'$ be the map induced by $f$ on $\operatorname{top}_y$. By the $\Lambda$-linearity of $f$, the square
\[
\begin{tikzcd}[column sep=large,row sep=large]
 W\oplus W\arrow[r,"d_{F_1(T)}"]\arrow[d,"h\oplus h"']
 &e_xF_1(T)\arrow[d,"f"]\\
 W'\oplus W'\arrow[r,"d_{F_1(U)}"']&e_xF_1(U)
\end{tikzcd}
\]
commutes. For each $v\in V$, we have $d_{F_1(T)}(a(v),b(v))=0$ by
\eqref{eq:recover-relation} for $T$, so the commutativity of the diagram gives $d_{F_1(U)}(ha(v),hb(v))=0$.
By \eqref{eq:recover-relation} for $U$ and the condition
$\ker a'\cap\ker b'=0$, there is a unique linear map $g:V\to V'$ satisfying $ha=a'g$ and $hb=b'g$.
We set $\rho_{T,U}(f)=(g,h)$. If $\rho_{U,T''}(f')=(g',h')$ for a composable morphism $f':F_1(U)\to F_1(T'')$,
then $f'f$ induces $h'h$ on $\operatorname{top}_y$, and $g'g$ satisfies the two defining equations
for $h'h$, so uniqueness yields
$\rho_{T,T''}(f'f)=(g'g,h'h)=\rho_{U,T''}(f')\rho_{T,U}(f)$. The same uniqueness argument
shows that $\rho_{T,U}$ is $k$-linear, and since
$F_1(g,h)$ is induced by $1\otimes h$, we have $\rho_{T,U}(F_1(g,h))=(g,h)$.
Finally, since $F_1(T)$ is generated by $e_yF_1(T)$, we have $\rho_{T,U}(f)=0$ if and only if $h=0$, which holds if and only if $f(F_1(T))\subseteq JF_1(U)$.
Every such homomorphism satisfies $f(J^iF_1(T))\subseteq
J^{i+1}F_1(U)$ for all $i\ge0$, so any $n$ of them compose to zero when $J^n=0$.
\end{proof}

\begin{example}\label{ex:not-full}
The embeddings of Theorem~\ref{thm:embedding} need not be full, and the
retractions of Proposition~\ref{prop:quotient-control} can have nonzero
kernels. Let
$\Lambda=k\langle \xi,\eta\rangle/(\xi^2,\xi\eta,\eta\xi,\eta^2)$
and take $s=\xi$, $t=\eta$ in Setting~\ref{set:quotient}.
This is a basic minimal representation-infinite algebra: its quotients
$\Lambda/(\xi-\lambda\eta)$, $\lambda\in k$, are pairwise non-isomorphic
indecomposable $\Lambda$-modules, whereas every proper nonzero algebra
quotient is isomorphic to $k$ or $k[t]/(t^2)$ and is representation-finite.
Taking $A=\Lambda$ and the identity reduction in Theorem~\ref{thm:embedding}
gives $F=F_1$. We have $F_1(S_-)=\Lambda$, and the retraction
$\End_\Lambda(\Lambda)\to\End_{K_2}(S_-)=k$ is the canonical projection
$\Lambda^{\mathrm{op}}\to k$, whose kernel is nonzero. Thus $F_1$ is not full.
\end{example}

\subsection{Vertex evaluations and block matrices}
We write $L=L_2$ or $L=L_3$ for the functors of
Constructions~\ref{con:cover} and~\ref{con:walk}, respectively, and
retain this notation for the rest of this section.
In the covering and walk constructions, we recover an endomorphism of $L(T)$ from the diagonal blocks of
an endomorphism of its push-down. The next two lemmas establish the
injectivity of vertex evaluations on $\End(L(T))$ and a nilpotence bound for the kernel
of taking the diagonal part. Finite direct sums in $\C$ will then give
retractions on all Hom spaces in Section~\ref{sec:uniform}.

\begin{lemma}\label{lem:vertex-evaluations}
Let $L$ be either functor just described and put
$\Sigma=\supp L(S_+)\cup\supp L(S_-)$. For nonzero $T\in\C$,
we have $\supp L(T)=\Sigma$, and evaluation
$\End(L(T))\to\End_k(L(T)(x))$ is injective at every $x\in\Sigma$;
here endomorphisms are taken over $R$ for the functor of
Construction~\ref{con:cover} and over $D_w$ for that of
Construction~\ref{con:walk}.
\end{lemma}
\begin{proof}
Write $T=(V,W;a,b)$. Applying the exact functor $L$ to
\eqref{eq:canonical-sequence} yields a short exact sequence. Since $L(S_+)$ and $L(S_-)$ have disjoint
supports, this sequence identifies $L(T)(x)$ naturally with $L(S_+)(x)\otimes V$ for
$x\in\supp L(S_+)$ and with $L(S_-)(x)\otimes W$ for $x\in\supp L(S_-)$.
Since $0\ne T\in\C$, both $V$ and $W$ are nonzero, so $\supp L(T)=\Sigma$.
Because $L$ is fully faithful, every endomorphism of $L(T)$ has the form $L(g,h)$ for some $(g,h)\in\End_{K_2}(T)$, and its component at $x\in\Sigma$
is $1\otimes g$ or $1\otimes h$. If this component vanishes, then $g=0$ or $h=0$,
and hence $(g,h)=0$ by Lemma~\ref{lem:Kronecker}(1).
\end{proof}

After summing vertex spaces over each fibre, morphisms between push-downs become block matrices.
Fern\'andez-Alonso, Gavito and P\'erez-Terrazas \cite{FernandezAlonso}
recover a morphism from a diagonal block of such a matrix in their
construction for truncated polynomial algebras; in
Propositions~\ref{prop:cover-end} and~\ref{prop:walk-end} below, we recover
endomorphisms from the diagonal part of such block matrices.
The next lemma uses injectivity of the vertex evaluations to establish a
nilpotence bound depending only on the maximum fibre size. In our
applications, $\Omega$ is the support of a module on the covering or the
vertex set of $D_w$, the partition $\mathcal P$ consists of the fibres over the vertices of $A$, and $\mathcal E$ is
the endomorphism algebra of the push-down.

\begin{lemma}
\label{lem:zero-column}
Let $\Omega$ be a finite nonempty set with a partition $\mathcal P$, and let
$M_i\ne0$ be vector spaces for $i\in\Omega$. Let $\mathcal E$ be a unital subalgebra
of $\End_k(\bigoplus_{i\in\Omega}M_i)$ preserving each subspace
$\bigoplus_{i\in\Omega'}M_i$, $\Omega'\in\mathcal P$. For a matrix
$f=(f_{j,i})$ with blocks $f_{j,i}:M_i\to M_j$, so that columns are indexed
by sources, let
$\delta(f)$ be the diagonal part of $f$, retaining only the blocks
$f_{i,i}$.
Assume:
\begin{enumerate}
\item $\delta(\mathcal E)\subseteq \mathcal E$;
\item for $\mathcal E_0=\delta(\mathcal E)$, every coordinate map
$\mathcal E_0\to\End_k(M_i)$ is injective;
\item $N=\ker(\delta|_{\mathcal E})$ is spanned by matrices $Z\in N$ such that
\begin{enumerate}[label={\upshape(\alph*)}]
\item each column of $Z$ has at most one nonzero block;
\item at least one column of $Z$ is zero.
\end{enumerate}
\end{enumerate}
Then $\delta:\mathcal E\to \mathcal E_0$ is a split algebra surjection, $N$ is a two-sided
ideal, and
\[
 N^m=0,\qquad m=\max_{\Omega'\in\mathcal P}|\Omega'|.
\]
\end{lemma}
\begin{proof}
By hypothesis (1), the diagonal matrices in $\mathcal E$ form a unital subalgebra $\mathcal E_0=\delta(\mathcal E)$, and $\mathcal E=\mathcal E_0\oplus N$ as vector spaces. Let $\mathcal Z$ be the set of matrices $Z\in N$ satisfying conditions (a) and (b) of hypothesis (3), so that $\mathcal Z$ spans $N$. Any product of matrices in $\mathcal Z$ again
has at most one nonzero block in each column, and a nonempty product
$Z$ of matrices in $\mathcal Z$ inherits a zero column from its rightmost factor. If column $i$ of $Z$ is zero,
then $Z_{i,i}=0$; since the coordinate map at $i$ is injective
on $\mathcal E_0$ by hypothesis (2), it follows that $\delta(Z)=0$.
In particular, $N^2\subseteq N$.
For $r\in \mathcal E_0$ and $f\in \mathcal E$, diagonal multiplication satisfies
$\delta(rf)=r\delta(f)$ and $\delta(fr)=\delta(f)r$,
so $\mathcal E_0N+N\mathcal E_0\subseteq N$. Thus
$N$ is a two-sided ideal of $\mathcal E$, and $\delta:\mathcal E\to \mathcal E_0$ is an algebra surjection split by the inclusion $\mathcal E_0\hookrightarrow \mathcal E$.

Now suppose that $Z_m\cdots Z_1$ is a product of $m$ matrices in $\mathcal Z$
having a nonzero block in column $i_0$. By condition (a), each factor $Z_t$
has at most one nonzero block in each column, so there is a unique sequence
of indices $i_0,i_1,\ldots,i_m\in\Omega$ such that $(Z_t)_{i_t,i_{t-1}}\ne0$
for all $1\le t\le m$, and the nonzero block of $Z_m\cdots Z_1$ in column
$i_0$ is the composite $(Z_m)_{i_m,i_{m-1}}\cdots(Z_1)_{i_1,i_0}$.
Since every factor preserves the subspaces indexed by the parts of $\mathcal P$,
all $m+1$ indices $i_0,\ldots,i_m$ belong to a single part of $\mathcal P$, which has at most $m$
elements. By the pigeonhole principle, there exist $0\le p<q\le m$ with $i_p=i_q$.
For the subproduct $Z'=Z_q\cdots Z_{p+1}$, column $i_p$ likewise has the
single nonzero block
\[
 Z'_{i_q,i_p}=(Z_q)_{i_q,i_{q-1}}\cdots(Z_{p+1})_{i_{p+1},i_p},
\]
which must be nonzero because the full composite
$(Z_m)_{i_m,i_{m-1}}\cdots(Z_1)_{i_1,i_0}$ is nonzero.
Since $i_p=i_q$, this is a nonzero diagonal block of $Z'$,
contradicting $\delta(Z')=0$ from the first paragraph.
Since $\mathcal Z$ spans $N$, we conclude that $N^m=0$.
\end{proof}

The nilpotence argument in Lemma~\ref{lem:zero-column} uses the same idea
as the proof of
\cite[Theorem~3.1]{GordonGreen}, where Gordon and Green show that a
$\mathbb Z$-graded Artin algebra whose subalgebra of degree zero is local
is itself local: there a repeated partial sum of degrees in a long product
produces a subproduct of degree zero, just as a repeated index produces a diagonal block
here.

\subsection{Endomorphisms of push-downs along a covering}
Let $\pi:R\to\mathcal A$ be a Galois covering with group $G$ as in
Definition~\ref{def:covering}, with push-down $\pi_*$ as in
Construction~\ref{con:pushdown}.
For objects $x,y$ in the same fibre of $\pi$, there is a
unique $\gamma\in G$ with $y=\gamma x$, and we group the blocks of an endomorphism according to this
group element, following the decomposition of morphisms between push-downs
in Gabriel \cite[\S3.2 and proof of Theorem~3.6(c)]{GabrielUniversal}.
The covering groups in Theorem~\ref{thm:finite-input}(2) are free and therefore
torsion-free, as required by the next proposition.

\begin{proposition}\label{prop:cover-end}
Suppose $G$ is torsion-free. Let $X$ be an $R$-module
with finite nonempty support $\Sigma$, and suppose every evaluation
$\End_R(X)\to\End_k(X(x))$, $x\in\Sigma$, is injective.
The map sending $f$ to the family of blocks $f_{x,x}:X(x)\to X(x)$
is a split algebra surjection
\[
 \End_A(\pi_*X)\longrightarrow\End_R(X)
\]
whose kernel has $m$-th power zero, where
$m=\max_q|\Sigma\cap\pi^{-1}(q)|$.
\end{proposition}
\begin{proof}
We decompose endomorphisms according to $G$, and then apply
Lemma~\ref{lem:zero-column} to the components indexed by $\gamma\ne1$.

\noindent\textbf{Step 1: Components of an endomorphism.}
Every $f\in\End_A(\pi_*X)$ preserves the fibre spaces
$(\pi_*X)(q)$, since it commutes with the vertex idempotents of $A$.
For $\gamma\in G$, let $f^{[\gamma]}$ be the matrix obtained by retaining only the blocks $f_{\gamma x,x}$. We show that each $f^{[\gamma]}$ belongs to $\End_A(\pi_*X)$.
Let $\psi:x\to y$ be a morphism in $R$. For each $\eta\in G$, the
decomposition of $\pi(\psi)$ with source $\eta x$ in
Construction~\ref{con:pushdown} has the single term $\eta(\psi):\eta x\to \eta y$, by
uniqueness in \eqref{eq:covering-lifts}; thus $(\pi_*X)(\pi(\psi))$ has
blocks $X(\eta(\psi)):X(\eta x)\to X(\eta y)$ for $\eta\in G$, and zero blocks elsewhere.
Comparing the blocks from $X(x)$ to $X(\gamma y)$ in
$f\circ(\pi_*X)(\pi(\psi))=(\pi_*X)(\pi(\psi))\circ f$, and using the fact that $G$ acts freely on
objects so that only $\eta=1$ contributes on the left and only $\eta=\gamma$ on the right,
we obtain
\begin{equation}\label{eq:cover-block}
 f_{\gamma y,y}X(\psi)=X(\gamma(\psi))f_{\gamma x,x}.
\end{equation}
Applying \eqref{eq:cover-block} to every $G$-translate of $\psi$ shows
that $f^{[\gamma]}$ commutes with $(\pi_*X)(\pi(\psi))$. Since the images of morphisms
of $R$ span the morphism spaces of $\mathcal A$, it follows that $f^{[\gamma]}\in\End_A(\pi_*X)$.
Because $\Sigma$ is finite, $f^{[\gamma]}$ is nonzero for only finitely many $\gamma\in G$.
For $\gamma=1$, equation~\eqref{eq:cover-block} is precisely the condition for an
endomorphism of $X$, so the diagonal part $f^{[1]}$ is the image of a unique element of
$\End_R(X)$ under $\pi_*$.

\noindent\textbf{Step 2: Nilpotence of the kernel.}
For $\gamma\ne1$, the matrix $f^{[\gamma]}$ has at most one nonzero block in
each row and column and has zero diagonal blocks. Moreover, $\gamma\Sigma\ne\Sigma$:
otherwise a positive power of $\gamma$ would fix an element of the finite
nonempty set $\Sigma$, contradicting the freeness of the action and the
torsion-freeness of $G$. Hence there exists $x\in\Sigma$ with $\gamma x\notin\Sigma$,
so the column of $f^{[\gamma]}$ indexed by $x$ is zero.
We now apply Lemma~\ref{lem:zero-column} to $\mathcal E=\End_A(\pi_*X)$, partitioning
$\Sigma$ by the fibres of $\pi$. Hypothesis (1) holds because
$\delta(f)=f^{[1]}\in \mathcal E$; hypothesis (2) is the assumed
injectivity of the vertex evaluations; and hypothesis (3) holds because
$f-\delta(f)=\sum_{\gamma\ne1}f^{[\gamma]}$. Lemma~\ref{lem:zero-column} therefore gives the result.
\end{proof}

\subsection{Endomorphisms of push-downs along a closed walk}
\label{sec:matching}
We work in Setting~\ref{set:walk} with the push-down functor $p_*$ of
Construction~\ref{con:walk}.
For an endomorphism of $p_*X$, the matrix blocks are indexed by pairs of
vertices lying in the same fibre of $p$, and the morphism relations link blocks
whose indices are connected by arrows of $D_w$ carrying the same label. We encode these index pairs and
relations in the matching quiver $H$ defined below, following the idea used by Crawley-Boevey
\cite[Theorem in \S2]{CrawleyBoeveyMaps} to describe homomorphisms between
tree modules and by Krause \cite[Theorem, p.~191]{KrauseTreeBand} between
tree and band modules via pairs of label-preserving maps from a common
quiver. Here the vertex spaces may be infinite-dimensional and the arrow maps $a,b$ out of $u$ may be noninvertible. We give the block argument in
this setting to obtain a nilpotence bound independent of their dimensions.

\begin{definition}\label{def:matching-quiver}
The \emph{matching quiver} $H$ has vertices $(j,i)$ with $p(j)=p(i)$,
ordered as target and source of a matrix block.
For each pair of arrows $\alpha:j\to j'$ and $\beta:i\to i'$ of $D_w$
with $p(\alpha)=p(\beta)$, it has an arrow
\[
 (\alpha,\beta):(j,i)\longrightarrow(j',i').
\]
Let $\Delta$ be the subquiver consisting of the diagonal vertices
$(i,i)$ and the arrows $(e,e)$, where $e$ runs over the arrows of $D_w$.
The coordinate projections $H\to D_w$ send $(j,i)$ to $j$ and $i$,
respectively, and $(\alpha,\beta)$ to $\alpha$ and $\beta$.
By a \emph{component} of $H$ we mean a connected component of its underlying unoriented graph.
\end{definition}

For an endomorphism $f:p_*X\to p_*X$, a vertex $(j,i)$ of $H$ indexes the
block $f_{j,i}:X(i)\to X(j)$, and
a matching arrow $(\alpha,\beta):(j,i)\to(j',i')$ encodes the relation
\begin{equation}\label{eq:walk-block}
 f_{j',i'}X(\beta)=X(\alpha)f_{j,i},
\end{equation}
which we derive in the proof of Proposition~\ref{prop:walk-end}.
For instance, for the walk $c(c')^{-1}$ on two distinct loops at a single
vertex, $D_w$ is the Kronecker quiver, the diagonal vertices form a
Kronecker component $\Delta$ of $H$, and the two off-diagonal vertices are isolated.
For the walk $w=\ell_1\ell_2\ell_1\ell_3^{-1}$ on three distinct loops at a single vertex, the
cycle $D_w$ and one off-diagonal component of $H$ are
\[
\begin{tikzcd}[column sep=large,row sep=large]
 d_0\arrow[r,"\ell_1"]\arrow[d,"\ell_3"']&d_1\arrow[d,"\ell_2"]\\
 d_3&d_2\arrow[l,"\ell_1"]
\end{tikzcd}
\qquad
\begin{tikzcd}[column sep=large]
 (d_2,d_0)\arrow[r]&(d_3,d_1).
\end{tikzcd}
\]
The two coordinate projections of this component pick out the two distinct arrows of $D_w$
labelled $\ell_1$, and each projection misses two vertices of $D_w$; accordingly, any block matrix
whose blocks vanish off the vertices of this component has zero columns at the vertices
missed by its second projection.
The next lemma establishes this property for all off-diagonal components of $H$, verifying hypothesis (3) of
Lemma~\ref{lem:zero-column} for Proposition~\ref{prop:walk-end}.

\begin{lemma}\label{lem:matching-components}
For a primitive cyclically reduced closed walk whose cycle $D_w$ has
a source, $\Delta$ is a whole component of $H$. Both coordinate
projections are injective on the vertices of each component. On a component
other than $\Delta$, neither projection is surjective onto $(D_w)_0$.
\end{lemma}
\begin{proof}
We first prove that the coordinate projections are injective on each component and that $\Delta$ is a component, then show that any component
with a surjective projection induces an automorphism of the cycle $D_w$, and
finally use the primitivity of $w$ to rule out such components other than $\Delta$.

\noindent\textbf{Step 1: Injectivity and the diagonal component.}
Put $r=|(D_w)_0|$. At any vertex $(j,i)$ of $H$, an arrow of $D_w$ incident to $j$
determines at most one arrow incident to $i$ with the same label and the same
orientation (incoming or outgoing), and conversely, by
\eqref{eq:local-injectivity}. Thus both coordinate projections $H\to D_w$ are
injective on the incident edges at each vertex, and every vertex of $H$ has degree at most two. Any path in $H$
that never traverses the same edge twice in succession therefore projects, in each coordinate, to a walk in the cycle $D_w$ that never
backtracks along the edge just traversed and hence moves in one fixed
direction around the cycle. If a shortest path in $H$ between two vertices had at least $r$
edges, its first $r$ edges would project to one complete turn in each
coordinate, returning to the starting vertex of $H$ and contradicting minimality; thus every shortest path in $H$ has fewer than $r$ edges.

Now suppose that two vertices in the same component of $H$ share a coordinate. The projection of a shortest path between them to that coordinate moves in one direction around $D_w$ and begins and ends at the same vertex, so its length is a multiple of $r$; since its length is less than $r$, it must be zero. This proves that both coordinate projections are injective on the vertices of each component.
Finally, at a diagonal vertex $(i,i)$, the two arrows in any incident edge of $H$
have the same label and the same source or target $i$, so they coincide by \eqref{eq:local-injectivity}; thus every edge of $H$ incident to $\Delta$ lies in $\Delta$, and $\Delta$ is a whole component.

\noindent\textbf{Step 2: A surjective projection gives a cycle automorphism.}
When $r=2$, since $D_w$ has a source $u$, both arrows of $D_w$ go from $u$ to the other vertex, so every off-diagonal
vertex of $H$ is isolated and cannot project onto $(D_w)_0$. Assume now that $r\ge3$, and let $C$ be a component of $H$ on which one of the coordinate projections is surjective. By Step~1, $C$ has at most $r$ vertices, so it has exactly $r$ vertices and both coordinate projections are bijective on the vertices of $C$.
Since $C$ is connected and has degree at most two, it is a path or a
cycle. It therefore contains a spanning path of length $r-1$, whose two projections are spanning paths
of the cycle $D_w$. Let $\sigma$ be the permutation of $(D_w)_0$ defined by requiring
$(\sigma(i),i)$ to lie in $C$ for every $i\in(D_w)_0$. Then $\sigma$ maps the second
projection of the spanning path edge by edge onto the first projection,
and therefore maps the two endpoints of the second projection to those of the first.
Since each projection omits a single edge of $D_w$, namely the edge joining
its two endpoints, $\sigma$ also maps the omitted edge of the second
projection to that of the first, and hence extends to an
automorphism of the underlying unoriented cycle of $D_w$.

\noindent\textbf{Step 3: Excluding a component other than $\Delta$.}
If $\sigma$ fixes a vertex $i$, then $(i,i)\in C$, so $C$ meets $\Delta$ and therefore equals $\Delta$ by Step~1.
Suppose for a contradiction that $C\ne\Delta$, so that $\sigma$ has no fixed vertex. Then $\sigma$ is either a nonidentity rotation or a reflection without fixed
vertices. A reflection without fixed vertices requires $r$ to be even;
it maps two opposite edges of the cycle to themselves while exchanging
their endpoints. Since the spanning path of $C$ omits only one edge, at least
one of these two edges, say $e:i\to j$, belongs to the second projection
of the spanning path. Here $\sigma(i)=j$ and $\sigma(j)=i$. By the
definition of a matching arrow, the edge of $C$ above $e$ would pair it
with an arrow $\sigma(i)\to\sigma(j)$, that is, an arrow $j\to i$. Since $r\ge3$, the cycle $D_w$ has only one arrow between $i$ and
$j$, so this is impossible.
Thus $\sigma$ is a nonidentity rotation.

Fix a cyclic orientation of $D_w$, which $\sigma$ preserves, and assign to each edge of $D_w$
its arrow label if its orientation agrees with that cyclic orientation, and the formal
inverse label otherwise. Along the spanning path of $C$, paired arrows have the same label and
orientation by Definition~\ref{def:matching-quiver}, so $\sigma$
preserves their signed labels. Consequently, the two projections of the spanning
path carry the same multiset of signed labels. Since adding the single omitted edge to
either projection yields all $r$ edges of $D_w$, the two omitted edges
must also carry the same signed label; since they correspond under $\sigma$ by Step~2,
the rotation $\sigma$ preserves the signed labels of all $r$ edges of $D_w$. If
$\sigma$ shifts vertices along the cycle by $\tau$ steps ($1\le\tau<r$), then the cyclic word of signed labels is
periodic with period $\ell=\gcd(r,\tau)<r$, and is therefore the
$(r/\ell)$-th power of its first $\ell$ letters. That shorter walk is
closed in $Q$, since the equal signed arrows at positions $1$ and $\ell+1$ have
the same initial vertex. Read from $d_0$, the full signed word is either $w$ or $w^{-1}$
according to the chosen cyclic orientation, so $w$ is a proper
power of a shorter closed walk, contradicting the primitivity of $w$.
\end{proof}

We now apply Lemma~\ref{lem:matching-components} to the endomorphism algebra
of a push-down $p_*X$.

\begin{proposition}\label{prop:walk-end}
Let $X$ be a representation of $D_w$ with nonzero
space at every vertex. Suppose $p_*X$ is an $A$-module and every
evaluation $\End_{D_w}(X)\to\End_k(X(i))$ is injective.
The map sending $f$ to the family of blocks $f_{i,i}:X(i)\to X(i)$
is a split algebra surjection
\[
 \End_A(p_*X)\longrightarrow\End_{D_w}(X)
\]
whose kernel has $m$-th power zero, where $m=\max_q|p^{-1}(q)|$.
\end{proposition}
\begin{proof}
We separate the morphism equations according to the components of $H$, and then apply
Lemma~\ref{lem:zero-column} to the components other than $\Delta$.

\noindent\textbf{Step 1: Decomposition by components of $H$.}
Since the relations for a morphism of $A$-modules are the same as those for
the underlying representations of $Q$, we may compute with the arrows of
$Q$. Let $f$ be an endomorphism of $p_*X$ and let $c:q\to q'$ be an arrow of $Q$.
For vertices $i\in p^{-1}(q)$ and $j'\in p^{-1}(q')$, the block from
$X(i)$ to $X(j')$ of the commutativity relation
$f\circ(p_*X)(c)=(p_*X)(c)\circ f$ has the form \eqref{eq:walk-block},
where $\beta:i\to i'$ is the unique arrow of $D_w$ labelled $c$ starting at $i$
and $\alpha:j\to j'$ is the unique arrow of $D_w$ labelled $c$ ending at $j'$, whenever
such arrows exist (with the corresponding term omitted if no such arrow exists). When both $\alpha$ and $\beta$ exist, the two blocks in \eqref{eq:walk-block} are indexed by
the endpoints $(j,i)$ and $(j',i')$ of the matching arrow $(\alpha,\beta)$ in $H$. When only $\beta$
exists, the relation reduces to $f_{j',i'}X(\beta)=0$, and here $j'\ne i'$ because
otherwise $\beta$ itself would be an arrow labelled $c$ ending at $j'$.
Likewise, when only $\alpha$ exists, the relation reduces to
$X(\alpha)f_{j,i}=0$ with $j\ne i$, since otherwise $\alpha$ would be an
arrow labelled $c$ starting at $i$. In both of these cases, the relation constrains
a single off-diagonal block. Consequently, every commutativity
relation involves only blocks indexed by vertices in a single component of
$H$, so retaining the blocks supported on any single component of $H$ yields an endomorphism of $p_*X$. In particular, the part supported on $\Delta$
is precisely the image of $\End_{D_w}(X)$ under $p_*$.

\noindent\textbf{Step 2: Nilpotence of the kernel.}
We apply Lemma~\ref{lem:zero-column} to $\mathcal E=\End_A(p_*X)$, partitioning the
vertices of $D_w$ by the fibres of $p$; $\mathcal E$ preserves the fibre
subspaces $(p_*X)(q)$. Hypothesis (1) holds because the diagonal part
of $f$ is the part supported on $\Delta$, which belongs to $\mathcal E$ by Step~1;
hypothesis (2) is the assumed injectivity of the vertex
evaluations; and hypothesis (3) holds because $\ker(\delta|_{\mathcal E})$
is spanned by the matrices obtained by restricting endomorphisms to
single components of $H$ other than $\Delta$. Each such matrix has
at most one nonzero block in each row and column because both coordinate
projections are injective on that component; it has a zero column because
the second projection is not surjective (Lemma~\ref{lem:matching-components}).
\end{proof}

\subsection{A uniform nilpotence bound}
\label{sec:uniform}
Since $\C$ is closed under finite direct sums
(Lemma~\ref{lem:Kronecker}(2)), Lemma~\ref{lem:vertex-evaluations} and
Propositions~\ref{prop:cover-end} and~\ref{prop:walk-end} apply to finite
direct sums in $\C$. Passing to matrix entries between direct summands
then yields retractions on all Hom spaces in $\C$ whose nilpotence bound
depends only on the fixed supports of $L(S_+)$ and $L(S_-)$, and is therefore
independent of the modules as required by \eqref{eq:nilpotence-bound}.

\begin{proposition}\label{prop:finite-control}
For each available alternative $\nu$ of Theorem~\ref{thm:finite-input},
the functor $F_\nu$ of \eqref{eq:chosen-functors} admits retractions
on $\C$ with nilpotence bound
\[
 n=\begin{cases}
 \min\{j\ge1:(\rad A)^j=0\},&\text{in alternative (1)},\\
 \max_q|\Sigma\cap\pi^{-1}(q)|,&\text{in alternative (2)},\\
 \max_q|p^{-1}(q)|,&\text{in alternative (3)},
 \end{cases}
\]
where $\Sigma=\supp L_2(S_+)\cup\supp L_2(S_-)$ in (2).
\end{proposition}
\begin{proof}
For alternative (1) of Theorem~\ref{thm:finite-input}, this follows from Proposition~\ref{prop:quotient-control}. For alternatives (2) and (3), fix $\nu\in\{2,3\}$; we first obtain
retractions on endomorphism algebras and then use finite direct sums to
extend them to all Hom spaces.

\noindent\textbf{Step 1: Endomorphism algebras.}
Put
$\Sigma=\supp L(S_+)\cup\supp L(S_-)$, and let $T\in\C$ be nonzero. In alternative (2),
$G$ is torsion-free and $L(T)$ has support $\Sigma$ with injective vertex evaluations by
Lemma~\ref{lem:vertex-evaluations}, so Proposition~\ref{prop:cover-end} applies to $L(T)$. In alternative (3),
$L_3(T)$ has full support $\Sigma=(D_w)_0$ and
injective vertex evaluations by Lemma~\ref{lem:vertex-evaluations}, and
$p_*L_3(T)$ is an $A$-module by Theorem~\ref{thm:finite-input}(3) and Lemma~\ref{lem:F3-properties},
so Proposition~\ref{prop:walk-end} applies.
In either case, taking the diagonal part and using the full faithfulness of $L$ gives a split algebra surjection
\begin{equation}\label{eq:finite-end-retraction}
 \End_A(F_\nu(T))\longrightarrow\End(L(T))
 \cong\End_{K_2}(T)
\end{equation}
whose kernel has $m$-th power zero, where
$m$ is the maximum fibre size in the fixed support $\Sigma$.

\noindent\textbf{Step 2: Hom spaces and the nilpotence bound.}
Let $T,U\in\C$.
If $T=0$ or $U=0$, then both Hom spaces vanish and we set
$\rho_{T,U}=0$.
Assume now that $T$ and $U$ are nonzero. Given $f:F_\nu(T)\to F_\nu(U)$, we
regard $f$ as the endomorphism of
$F_\nu(T\oplus U)=F_\nu(T)\oplus F_\nu(U)$
whose $(U,T)$-component is $f$ and whose other three components are zero; we then apply
\eqref{eq:finite-end-retraction} for $T\oplus U$ to this endomorphism and define
$\rho_{T,U}(f):T\to U$ to be the $(U,T)$-component of the resulting
endomorphism of $T\oplus U$. Explicitly, at each vertex $i\in\Sigma$, the $(U,T)$-component of the $(i,i)$-block retained by
\eqref{eq:finite-end-retraction} is
simply the $(i,i)$-block $f_{i,i}:L(T)(i)\to L(U)(i)$ of $f$. Because the inclusions and
projections for a finite direct sum act separately at each vertex,
the resulting family of vertex maps $(f_{i,i})_{i\in\Sigma}:L(T)\to L(U)$ is independent of the ambient finite direct sum in $\C$ used to compute it,
and $\rho_{T,U}(f)$ is its unique preimage under the fully faithful functor $L$.
In particular, $\rho_{T,U}$ is $k$-linear, and since the $(i,i)$-blocks of
$F_\nu(g,h)$ are the vertex maps of
$L(g,h)$, we have $\rho_{T,U}(F_\nu(g,h))=(g,h)$.
Moreover, for a third nonzero object $T''\in\C$ and a morphism $f':F_\nu(U)\to F_\nu(T'')$,
applying the algebra surjection \eqref{eq:finite-end-retraction} to $T\oplus U\oplus T''$
shows that $\rho_{T,T''}(f'f)=\rho_{U,T''}(f')\rho_{T,U}(f)$.
Finally, for $U=T$, this construction yields the endomorphism retraction
of Step~1, whose kernel has $m$-th power zero; hence these retractions have
nilpotence bound $n=m$.
\end{proof}

\section{The embedding theorem and large families}
\label{sec:completion}
Throughout this section, $k$ is algebraically closed and $A$ is an
arbitrary representation-infinite finite-dimensional $k$-algebra.
We finish the embedding theorem by cancelling the simple projective
summands and transferring the functors from the reduced algebra.
We then bound the nilpotence exponent and apply the embedding to
G\"obel--Simson families.

The cancellation uses the following fact about a module with local
endomorphism ring and nilpotent radical. Maps factoring through direct sums of copies of this module form a
nilpotent ideal on modules having no summand isomorphic to it.
\begin{lemma}\label{lem:fixed-summand}
Let $P$ be a nonzero $A$-module with local endomorphism ring and choose
$\ell\ge1$ with $(\rad\End_A(P))^\ell=0$.
On the full subcategory of $\Mod A$ consisting of modules with no direct
summand isomorphic to $P$, let $\mathcal J(X,Y)$ consist of maps
$X\to Y$ factoring through an arbitrary direct sum of copies of $P$.
Then $\mathcal J$ is an ideal and $\mathcal J^{\ell+1}=0$.
\end{lemma}
\begin{proof}
Closure under addition and composition with arbitrary morphisms is immediate
from the definition of $\mathcal J$. For any module $X$ with no direct summand isomorphic to $P$, every composite $P\to X\to P$ is a nonunit of
$\End_A(P)$, since a unit would split a copy of $P$ off $X$.
Since $\End_A(P)$ is local, every such composite belongs to $\rad\End_A(P)$.

Now consider $\ell+1$ composable morphisms in $\mathcal J$, written as
$X_{i-1}\xrightarrow{u_i}P^{(I_i)}\xrightarrow{v_i}X_i$ ($1\le i\le \ell+1$).
Their composite is
\[
 v_{\ell+1}(u_{\ell+1}v_\ell)\cdots(u_2v_1)u_1.
\]
For each $1\le i\le \ell$, put $w_i=u_{i+1}v_i:P^{(I_i)}\to P^{(I_{i+1})}$;
by the first paragraph, every matrix component
$(w_i)_{s,r}:P\to P$ ($r\in I_i$, $s\in I_{i+1}$) belongs to
$\rad\End_A(P)$. Fix $s_1\in I_1$ and an element $x\in P$ in the
$s_1$-th summand of $P^{(I_1)}$. Since each $w_i$ sends an element of a
summand of $P^{(I_i)}$ to an element of $P^{(I_{i+1})}$ with only finitely
many nonzero coordinates, the $s_{\ell+1}$-th coordinate of
$w_\ell\cdots w_1(x)$ in $P^{(I_{\ell+1})}$ is the finite sum
\[
 \sum_{s_2\in I_2,\ldots,s_\ell\in I_\ell}
 (w_\ell)_{s_{\ell+1},s_\ell}\cdots(w_1)_{s_2,s_1}(x)
\]
for every $s_{\ell+1}\in I_{\ell+1}$. Each product
$(w_\ell)_{s_{\ell+1},s_\ell}\cdots(w_1)_{s_2,s_1}$ lies in
$(\rad\End_A(P))^\ell=0$, so every term of the sum vanishes.
Since elements supported in single summands span $P^{(I_1)}$, we have
$w_\ell\cdots w_1=0$, and hence $v_{\ell+1}w_\ell\cdots w_1u_1=0$.
\end{proof}

We can now prove Theorem~A for every choice in
Theorem~\ref{thm:finite-input}.
\begin{theorem}\label{thm:embedding}
Let $A$ be a representation-infinite finite-dimensional algebra over an
algebraically closed field $k$. Let $A_0$, $B$ and $\Phi$ be as in
Lemma~\ref{lem:reduction}, and let
$F_B:\D\to\Mod B$ be any of the functors $F_\nu$ in
\eqref{eq:chosen-functors}, formed from data in
Theorem~\ref{thm:finite-input} for $B$.
Then $F=\Phi F_B:\D\to\Mod A$ is a
faithful representation embedding with the following properties.
\begin{enumerate}
\item $F$ sends finite-dimensional modules to finite-dimensional modules.
\item $F$ commutes with arbitrary direct sums and filtered colimits.
\item If $T\in\D$ has infinite dimension, then $\dm F(T)=\dm T$.
\end{enumerate}
Moreover, $F$ admits retractions on $\C$ with the same nilpotence bound
as $F_B$ in Proposition~\ref{prop:finite-control}.
\end{theorem}
\begin{proof}
Lemmas~\ref{lem:F1-properties}, \ref{lem:L2}, \ref{lem:F3-properties}
and \ref{lem:reduction} show that $F$ is exact and faithful.
We first check the three stated properties.

(1) The same lemmas show that $F_B$ and $\Phi$ preserve finite-dimensional modules.

(2) Those lemmas also show that $F_B$ and $\Phi$ commute with arbitrary
direct sums and filtered colimits.

(3) By Lemmas~\ref{lem:F1-properties}, \ref{lem:L2} and
\ref{lem:F3-properties}, $\dm F_B(T)=\dm T$ for infinite-dimensional
$T$. Lemma~\ref{lem:reduction} gives $\dm F(T)=\dm F_B(T)$.

For $T,U\in\C$, compose the inverse of the Hom isomorphism induced by
$\Phi$ with the retraction for $F_B$:
\[
 \Hom_A(F(T),F(U))\xrightarrow{\sim}
 \Hom_B(F_B(T),F_B(U))\xrightarrow{\rho^B_{T,U}}\Hom_{K_2}(T,U).
\]
These maps preserve composition and have the bound in
Proposition~\ref{prop:finite-control}. By
Lemma~\ref{lem:retraction-consequences}, $F|_\C$ preserves
indecomposability and reflects isomorphism classes.
Put $F_{\mathrm{fin}}=F|_{\D_{\mathrm{fin}}}$. The finite restriction
of $F_B$ is a representation embedding by Theorem~\ref{thm:finite-input}
and, in alternative (2), Lemma~\ref{lem:L2-identification}. Since $\Phi$
is exact and fully faithful by Lemma~\ref{lem:reduction},
$F_{\mathrm{fin}}$ is also a representation embedding. Thus $P=F(S_-)$
is finite-dimensional and indecomposable, so $\End_A(P)$ is local with
nilpotent radical. We complete the proof on $\D$ in three steps.

\noindent\textbf{Step 1: Excluding $P$ as a summand.}
Let $T_0\in\C$, and suppose that $P\to F(T_0)$ is a split monomorphism with retraction $p:F(T_0)\to P$. By part (2) and \eqref{eq:colimit-comparison}, $F(T_0)$ is the filtered
colimit of the modules $F(U)$ for finite-dimensional $U\subseteq T_0$. Since $P$
is finitely presented and the finite-dimensional submodules $T_1\subseteq T_0$ belonging to
$\C_{\mathrm{fin}}$ are cofinal
(Lemma~\ref{lem:Kronecker}(3)), the map $P\to F(T_0)$ factors through
$F(T_1)$ for some $T_1\in\C_{\mathrm{fin}}$. Composing the canonical map
$F(T_1)\to F(T_0)$ with $p$ shows that $P$ is isomorphic to a direct summand of
$F(T_1)$.
The indecomposable summands of $T_1$ belong to $\C_{\mathrm{fin}}$,
since they have neither simple summand, and their images under $F$
are indecomposable. The Krull--Schmidt theorem therefore gives an
indecomposable summand $T'$ of $T_1$ with $F(T')\cong P=F(S_-)$.
Since $F_{\mathrm{fin}}$ reflects isomorphism classes,
$T'\cong S_-$,
contradicting the fact that $T_1\in\C$ has no summand isomorphic to $S_-$.

\noindent\textbf{Step 2: Reflection of isomorphism classes.}
Let $T,U\in\D$ satisfy $F(T)\cong F(U)$. By
\eqref{eq:simple-splitting}, we may write $T=T_0\oplus S_-^{(I)}$ and
$U=U_0\oplus S_-^{(K)}$ with $T_0,U_0\in\C$. Put $X=F(T_0)$ and
$Y=F(U_0)$. Since $F$ commutes with direct sums, an
isomorphism $F(T)\xrightarrow{\sim}F(U)$ and its inverse take the block form
\[
 f=\begin{pmatrix}f_{11}&f_{12}\\f_{21}&f_{22}\end{pmatrix}:
 X\oplus P^{(I)}\longrightarrow Y\oplus P^{(K)},\qquad
 f^{-1}=\begin{pmatrix}g_{11}&g_{12}\\g_{21}&g_{22}\end{pmatrix},
\]
so that
\begin{equation}\label{eq:cancellation}
 g_{11}f_{11}=1_X-g_{12}f_{21},\qquad
 f_{11}g_{11}=1_Y-f_{12}g_{21}.
\end{equation}
The composites $g_{12}f_{21}$ and $f_{12}g_{21}$ factor through direct sums of
copies of $P$. Since $X$ and $Y$ have no direct
summand isomorphic to $P$ by Step~1, Lemma~\ref{lem:fixed-summand} implies that
$g_{12}f_{21}$ and $f_{12}g_{21}$ are nilpotent. Consequently, both
right sides of \eqref{eq:cancellation} are invertible (with inverses given
by finite geometric series), so $f_{11}:X\to Y$ is an isomorphism, and
Lemma~\ref{lem:retraction-consequences} gives $T_0\cong U_0$.
Moreover, multiplying $f$ on the left and right by invertible triangular matrices
yields
\[
 \begin{pmatrix}1&0\\-f_{21}f_{11}^{-1}&1\end{pmatrix}
 f\begin{pmatrix}1&-f_{11}^{-1}f_{12}\\0&1\end{pmatrix}
 =\begin{pmatrix}f_{11}&0\\0&f_{22}-f_{21}f_{11}^{-1}f_{12}\end{pmatrix},
\]
which shows that $P^{(I)}\cong P^{(K)}$. Taking dimensions over $k$ gives $(\dm P)|I|=(\dm P)|K|$;
since $0<\dm P<\infty$, we obtain $|I|=|K|$ and hence $T\cong U$.

\noindent\textbf{Step 3: Preservation of indecomposability.}
By
\eqref{eq:simple-splitting}, every indecomposable object of $\D$ either
belongs to $\C$ or is isomorphic to $S_-$.
In the first case, $F$ preserves indecomposability by Lemma~\ref{lem:retraction-consequences}; in the
second, $F(S_-)\cong P$ is indecomposable because $F_{\mathrm{fin}}$ preserves indecomposability.
\end{proof}

\begin{remark}\label{rem:weaker-finite-input}
In the proof of Theorem~\ref{thm:embedding}, the assumption that
$F_{\mathrm{fin}}$ is a representation embedding can be replaced by
the indecomposability of $F(S_-)$ and the condition
$F(S_-)\not\cong F(U)$ for every $U\in\C_{\mathrm{fin}}$, with the
other properties of $F$ and its retractions retained.
These are the two finite-dimensional facts used in the Lean proof of
Theorem~\ref{thm:embedding}. They also suffice for the proof given here:
the retractions imply that $F|_\C$ preserves indecomposability, so the
Krull--Schmidt argument in Step~1 would give $U\in\C_{\mathrm{fin}}$
with $F(U)\cong F(S_-)$, contrary to the stated condition. Step~3 uses
only the indecomposability of $F(S_-)$.
\end{remark}

The rank bound in \cite[Theorem~7]{Bongartz} will bound the nilpotence exponent.
We record it together with the choices of construction made in that proof.
Write
$\iota:\md k[t]\to\D_{\mathrm{fin}}$ for the functor
$\iota(U)=(U,U;t,1)$, acting by $(f,f)$ on a $k[t]$-module homomorphism $f$.
\begin{lemma}[{\cite[\S5.2, Theorem~7 and \S5.4]{Bongartz}}]
\label{lem:quantitative-input}
Let $B$ be a connected basic minimal representation-infinite algebra
over an algebraically closed field $k$, and put $d=\dm B$.
\begin{enumerate}
\item If the lattice of two-sided ideals of $B$ is not distributive,
then alternative (1) of Theorem~\ref{thm:finite-input} applies.
\item If this lattice is distributive, there is a choice of $F_B$ in
alternative (2) or (3) of Theorem~\ref{thm:finite-input} and a bimodule
${}_B\mathcal B_{k[t]}$, free as a right $k[t]$-module of rank
$1\le \mu\le\max(30,4d)$, such that
\[
 (F_B|_{\D_{\mathrm{fin}}})\circ\iota\cong \mathcal B\otimes_{k[t]}-
 \quad\text{as functors }\md k[t]\longrightarrow\md B.
\]
\end{enumerate}
\end{lemma}

Part (1) is the final paragraph of \cite[\S5.2]{Bongartz}, applying
\cite[Corollary~2]{Bongartz}. The tensor description and the upper bound
in part (2) are the assertion of \cite[Theorem~7]{Bongartz}.
The proof in \cite[\S5.4]{Bongartz} obtains this tensor functor from
one of the covering or walk constructions in
Theorem~\ref{thm:finite-input}. The lower bound $\mu\ge1$ follows because
the composite is a representation embedding. Evaluating the tensor
functor at $k[t]/(t)$ bounds the support sizes in
Proposition~\ref{prop:finite-control} and hence gives the next result.
\begin{proposition}\label{prop:global-refinement}
There is a choice of the embedding $F$ in Theorem~\ref{thm:embedding}
and an integer $1\le n\le\max(30,4\dm A)$ such that $F$ admits
retractions on $\C$ with nilpotence bound $n$.
\end{proposition}
\begin{proof}
Choose $A_0=eAe$ and $B$ as in Lemma~\ref{lem:reduction}, and put
$d=\dm B\le\dm A$. If the lattice of two-sided ideals of $B$ is not
distributive, Lemma~\ref{lem:quantitative-input}(1) gives alternative (1). Here the quotient construction is over $B$ itself, whose Loewy length
is at most $d$, so the assertion follows
from Proposition~\ref{prop:finite-control}.

If the lattice of two-sided ideals of $B$ is distributive, choose
$F_B$ and the bimodule $\mathcal B$, free of rank $\mu\le\max(30,4d)$ over $k[t]$, as
in Lemma~\ref{lem:quantitative-input}(2). For
$T_*=(k,k;0,1)=\iota(k[t]/(t))$, we have
$F_B(T_*)\cong \mathcal B\otimes_{k[t]}k[t]/(t)$ and hence $\dm F_B(T_*)=\mu$.
In alternative (2), every vertex of
$\Sigma=\supp M_-\cup\supp M_+$ contributes at least one dimension
to $F_B(T_*)$, so the maximum fibre size is at most
$|\Sigma|\le \mu$. In alternative (3), all $r$ vertex spaces of
$L_3(T_*)$ are one-dimensional, so $r=\mu$ and the maximum fibre size
is at most $\mu$. By Proposition~\ref{prop:finite-control} and Theorem~\ref{thm:embedding}, the resulting functor $F=\Phi F_B$ admits retractions on $\C$ with nilpotence bound $n\le\mu\le\max(30,4\dm A)$.
\end{proof}

We use the following specialization of the rigid families of
G\"obel and Simson. It records the cardinality and Hom spaces needed
for the application.
\begin{theorem}[{\cite[Theorem~1.2(i),(iii) and Theorem~3.3]{GoebelSimson}}]
\label{thm:rigid-family}
For every field $k$ and every infinite cardinal $\kappa$, there are
Kronecker modules $T_\beta=(V_\beta,W_\beta;a_\beta,b_\beta)$,
indexed by $\beta\subseteq\kappa$, such that
\[
 \dm V_\beta=\dm W_\beta=\kappa,
 \qquad a_\beta,b_\beta\text{ are injective},\qquad
 a_\beta(V_\beta)+b_\beta(V_\beta)=W_\beta,
\]
and
\[
 \Hom_{K_2}(T_\beta,T_\gamma)\cong
 \begin{cases}k,&\beta\subseteq\gamma,\\0,&\beta\nsubseteq\gamma,\end{cases}
 \qquad \End_{K_2}(T_\beta)\cong k.
\]
The Hom isomorphisms are $k$-linear and the endomorphism isomorphisms
are isomorphisms of $k$-algebras.
\end{theorem}
The displayed arrow conditions show that $T_\beta\in\C$.

Applying our chosen embedding gives the announced families.
\begin{corollary}\label{cor:radical}
Let $A$ be a representation-infinite finite-dimensional algebra over an
algebraically closed field $k$. There is an integer
$1\le n\le\max(30,4\dm A)$ such that, for every infinite cardinal
$\kappa$, there are $2^\kappa$ pairwise non-isomorphic indecomposable
$A$-modules $M$ satisfying
\[
 \dm M=\kappa,\qquad
 \End_A(M)/\rad\End_A(M)\cong k,\qquad
 (\rad\End_A(M))^n=0.
\]
If $\kappa\ge |k|$, this is the largest possible number of isomorphism
classes in dimension $\kappa$.
\end{corollary}
\begin{proof}
Fix $F$ and $n$ as in Proposition~\ref{prop:global-refinement}, and
put $M_\beta=F(T_\beta)$ for the family in
Theorem~\ref{thm:rigid-family}. Each $T_\beta$ is indecomposable
because $\End_{K_2}(T_\beta)\cong k$, and for $\beta\ne\gamma$ at
least one of $\Hom_{K_2}(T_\beta,T_\gamma)$ and $\Hom_{K_2}(T_\gamma,T_\beta)$ vanishes,
so the modules $T_\beta$ are pairwise non-isomorphic. Theorem~\ref{thm:embedding} therefore
gives $2^\kappa$ pairwise non-isomorphic indecomposable $A$-modules $M_\beta$ of
dimension $\kappa$.
For each $\beta$, the retraction $\rho_{T_\beta,T_\beta}$ is a surjective $k$-algebra homomorphism
\[
 \End_A(M_\beta)\longrightarrow\End_{K_2}(T_\beta)\cong k
\]
whose kernel $N_\beta$ satisfies $N_\beta^n=0$.
Every nilpotent ideal is contained in the Jacobson radical, and the quotient $k$
has zero radical, so
$N_\beta=\rad\End_A(M_\beta)$.

Finally, suppose that $\kappa\ge|k|$. Any $k$-vector space of dimension $\kappa$
has cardinality $\kappa$ and therefore admits at most
$\kappa^\kappa=2^\kappa$ linear endomorphisms. Since an $A$-module structure
on this space is determined by the actions of a finite basis of $A$, there are
at most $(2^\kappa)^{\dm A}=2^\kappa$ such module structures, and hence at most $2^\kappa$
isomorphism classes of $A$-modules of dimension $\kappa$.
\end{proof}

Note that in Corollary~\ref{cor:radical}, the conclusion
$\End_A(M)/\rad\End_A(M)\cong k$ cannot in general be replaced by $\End_A(M)\cong k$.
Indeed, over any finite-dimensional local algebra $A$ with $A/\rad A\cong k$ (for example, the representation-infinite local algebra $k\langle x,y\rangle/(x,y)^2$), every nonsimple indecomposable $A$-module $M$ satisfies
$(\rad A)M\ne0$, since otherwise $M$ would be semisimple. Because $\rad A$ is nilpotent, both the top $M/(\rad A)M$
and the socle of $(\rad A)M$ are nonzero; composing a nonzero projection
$p:M\to k$ factoring through $M/(\rad A)M$ with an injection $i:k\to M$ whose image lies in
$(\rad A)M$ gives a nonzero endomorphism $ip\in\End_A(M)$ with $(ip)^2=0$.

\appendix
\section{Corner quotients of a special biserial algebra}
\label{sec:simson-counterexample}
Throughout this appendix, $k$ is an arbitrary field and all algebras are
finite-dimensional $k$-algebras.
We examine the assertion about corner quotients used in Case~4 of the proof of
\cite[Theorem~3.1, pp.~128--129]{SimsonCorner}.
Proposition~3.7 of \cite{SimsonCorner}
asserts that every representation-infinite special biserial algebra $R$
over a field $k$ admits an idempotent $e$ and a surjective $k$-algebra homomorphism
$eRe\to S$, where $S$ is either
$k\langle x,y\rangle/(x^2,xy,yx,y^2)$ or the path algebra of an acyclic
orientation of a cycle with at least two vertices; here $eRe$ is the
\emph{corner} of $R$ determined by $e$. Case~4 applies this assertion
to an arbitrary representation-infinite special biserial algebra and then
transfers modules from $S$ to $R$ through $eRe$.
Example~\ref{ex:simson-corner} below shows that over any field $k$, the required corner quotient
need not exist. Because the algebra $R$ in Example~\ref{ex:simson-corner} has radical square
zero, Case~1 of the same proof already establishes the conclusion of \cite[Theorem~3.1]{SimsonCorner} for $R$; thus the example refutes Proposition~3.7, invoked in Case~4,
rather than the conclusion of \cite[Theorem~3.1]{SimsonCorner}.

Recall that a finite-dimensional algebra $kQ/I$ with $I$ admissible is
\emph{special biserial} if each vertex has at most two incoming and two
outgoing arrows and, for every arrow, at most one preceding arrow and
at most one following arrow have a product outside $I$
\cite[p.~124, (SP1)--(SP2)]{SimsonCorner}.
The following algebra satisfies these conditions over any field $k$ but has no corner
admitting either quotient above.

\begin{example}\label{ex:simson-corner}
Let $k$ be any field. Let $Q$ have vertices $1,2$, loops
$\alpha$ at $1$ and $\delta$ at $2$, and arrows
$\beta:1\to2$ and $\gamma:2\to1$. Put
\[
 R=kQ/(\text{all paths of length two}).
\]
Then $R$ is a six-dimensional representation-infinite special biserial
algebra. For every idempotent $e\in R$, the corner $eRe$ has no quotient
isomorphic to
\[
 k\langle x,y\rangle/(x^2,xy,yx,y^2)
\]
or to the path algebra of an acyclic orientation of a cycle with at
least two vertices (for two vertices, this is the Kronecker quiver).
\end{example}
\begin{proof}
Let $e_i$ be the trivial path at vertex $i$ for $i=1,2$; the elements $e_1,e_2,\alpha,\beta,\gamma,\delta$ form a basis of $R$.
Since each vertex has two incoming and two outgoing arrows and all products
of two arrows vanish in $R$, the algebra $R$ is special biserial.
For each integer $m\ge1$, let $C_m$ be the nilpotent Jordan block of
size $m$, and define a representation $M_m$ of $Q$ by placing $k^m\oplus k^m$ at each vertex and
setting
\[
 M_m(\alpha)=M_m(\delta)=M_m(\beta)=
 \begin{pmatrix}0&0\\I_m&0\end{pmatrix},
 \qquad
 M_m(\gamma)=\begin{pmatrix}0&0\\C_m&0\end{pmatrix}.
\]
Since all products of two arrow maps vanish, $M_m$ is an $R$-module.
Any endomorphism of $M_m$ commutes with the actions of $\alpha$ and $\delta$, so its vertex matrices take the
form
\[
 E_i=\begin{pmatrix}X_i&0\\Y_i&X_i\end{pmatrix}\quad(i=1,2).
\]
Commutativity with $\beta$ gives $X_1=X_2=:X$, and commutativity with $\gamma$ gives
$XC_m=C_mX$. Because $C_m$ has a cyclic vector, its commutant in $\operatorname{Mat}_m(k)$ is
$k[C_m]\cong k[t]/(t^m)$. The algebra homomorphism
$\End_R(M_m)\to k[C_m]$ sending $(E_1,E_2)$ to $X$ is surjective
with square-zero kernel; since $k[C_m]$ is local, $\End_R(M_m)$ is also
local, and hence $M_m$ is indecomposable. Since $\dm M_m=4m$ is
unbounded, $R$ is representation-infinite over every field $k$.

Put $J=\rad R$, so that $J^2=0$. Every idempotent $e\in R$ is conjugate to one
of $0,e_1,e_2,1$: writing $e=f+r$ with
$f\in\{0,e_1,e_2,1\}$ and $r\in J$, the relation $e^2=e$ gives
$r=fr(1-f)+(1-f)rf$, so $z=(1-f)rf-fr(1-f)$ satisfies $z^2=0$ and
\[
 (1+z)f(1-z)=f+r=e.
\]
Up to isomorphism, the nonzero corners of $R$ are therefore
$e_1Re_1\cong e_2Re_2\cong k[t]/(t^2)$ and $R$ itself. The two-dimensional corners $k[t]/(t^2)$ cannot map onto any algebra of dimension at least three. Under any algebra homomorphism from $R$ to the local three-dimensional
algebra $k\langle x,y\rangle/(x^2,xy,yx,y^2)$, the orthogonal idempotents $e_1,e_2$
must map to $1,0$ in some order; consequently, $\beta$, $\gamma$ and the loop at the
vertex whose idempotent maps to $0$ lie in the kernel, so the image is spanned by $1$
and the remaining loop and no such homomorphism can be surjective.

Finally, every quotient of a corner of $R$ has at most two isomorphism classes of simple
modules, which rules out path algebras of cycle quivers with three or more
vertices. The only remaining candidate is the
Kronecker algebra $kK_2$. Any algebra surjection $R\to kK_2$ would map
$e_1,e_2$ to the two primitive orthogonal idempotents of $kK_2$ in some order;
however, both $e_1Re_2$ and $e_2Re_1$ have dimension one, whereas the arrow space of $kK_2$ has dimension two, so no such surjection exists.
\end{proof}

\end{document}